\documentclass[leqno,11pt]{article}  
\usepackage{amsmath}
\usepackage{amssymb}

\usepackage[german,english]{babel}
\usepackage{amsthm}

\title{Sheaves of bounded $p$-adic logarithmic differential forms}

\author{\textsc{Elmar Grosse-Kl\"onne}}
\date{}

\theoremstyle{plain} 
\newtheorem{satz}{Theorem}[section]  
\newtheorem{kor}[satz]{Corollary}  
\newtheorem{lem}[satz]{Lemma}  
\newtheorem{pro}[satz]{Proposition}  
\newcommand{\ho}{\mbox{\rm Hom}}  
 
\newcommand{\spec}{\mbox{\rm Spec}}  
\newcommand{\proj}{\mbox{\rm Proj}}
\newcommand{\spf}{\mbox{\rm Spf}}  

\newcommand{\bi}{\mbox{\rm im}}  
\newcommand{\ke}{\mbox{\rm Ker}}  
\newcommand{\gr}{\mbox{\rm gr}}

\newcommand{\kara}{\mbox{\rm char}}  

\newcommand{\dlog}{\mbox{\rm dlog}}

\newcommand{\diag}{\mbox{\rm diag}}

\theoremstyle{remark}

\theoremstyle{definition}

\newcommand{\F}{\ensuremath{\mathbb{L}}}

\DeclareMathOperator{\Hom}{Hom}

\newcommand{\0}{\ensuremath{\overrightarrow{0}}}

\begin{document}
\maketitle
\footnote[0]{I wish to thank Ahmed Abbes, Adrian Iovita, Arthur Ogus, Peter Schneider, Ehud de Shalit and Matthias Strauch for their interest in this work. I am grateful to the referee for his very careful reading.}

\begin{abstract}

Let $K$ be a local field, $X$ the Drinfel'd symmetric space $X$ of dimension $d$ over $K$ and ${\mathfrak X}$ the natural formal ${\mathcal O}_K$-scheme underlying $X$; thus $G={\rm GL}\sb {d+1}(K)$ acts on $X$ and ${\mathfrak X}$. Given a $K$-rational $G$-representation $M$ we construct a $G$-equivariant subsheaf ${\mathcal M}^0_{{\mathcal O}_{\dot{K}}}$ of ${\mathcal O}_K$-lattices in the constant sheaf $M$ on ${\mathfrak X}$. We study the cohomology of sheaves of logarithmic differential forms on $X$ (or ${\mathfrak X}$) with coefficients in ${\mathcal M}^0_{{\mathcal O}_{\dot{K}}}$. In the second part we give general criteria for two conjectures of P. Schneider on $p$-adic Hodge decompositions of the cohomology of $p$-adic local systems on projective varieties uniformized by $X$. Applying the results of the first part we prove the conjectures in certain cases.  
\end{abstract}

%


\begin{center} {\bf Introduction} \end{center}

Let $p$ be prime number and $d\in\mathbb{N}$, let $K/{\mathbb Q}_p$ be a finite extension. In connection with the search for a Langlands type correspondence between suitable $p$-adically continuous representations of the group ${\rm GL}_{d+1}(K)$ on $p$-adic vector spaces on the one hand, and suitable $p$-adic Galois representations on the other hand, the $p$-adic cohomology (de Rham, crystalline, coherent, $p$-adic \'{e}tale) of Drinfel'd's symmetric space $X$ over $K$ and its projective quotients $X_{\Gamma}=\Gamma\backslash X$, with coefficients in rational representations $M$ of ${\rm GL}_{d+1}(K)$, has recently found increasing interest. We mention the first spectacular results due to Breuil \cite{br} who uses the cohomology of $X$ and $X_{\Gamma}$ with coefficients in $M={\rm Sym}^k({\mathbb Q}_p^2)$ (some $k\in{\mathbb N}$) to establish a partial correspondence in case $d=1$, $K={\mathbb Q}_p$, and the work of Schneider and Teitelbaum \cite{ast} where (for any $d$ and $K$) the ${\rm GL}_{d+1}(K)$-representation on the space $\Omega_X^{\bullet}(X)$ of top differential forms on $X$ is determined. Substantial as these works are, they call for generalizations. On the one hand one hopes to generalize the constructions from \cite{br} to cases where $d>1$. Since a decisive ingredient in \cite{br} is the work with $p$-adic {\it integral structures} in equivariant sheaf complexes on $X$, the investigation of such integral structures should be a starting point. On the other hand one hopes to generalize the analysis of \cite{ast} to more general equivariant vector bundles on $X$ (instead of the line bundle $\Omega_X^d$), e.g. to the vector bundles $\Omega_X^i$, or even $M\otimes\Omega_X^i$, for any $i$; this would be done best by finding and analysing equivariant subsheaves in $M\otimes\Omega_X^i$, e.g. those of exact, closed or logarithmic differential forms.
All this motivates the main objetive of the first part of this paper, the study of equivariant integral structures in the vector bundles $M\otimes\Omega_X^i$ and in their sub sheaves of logarithmic differential forms. The central question concerning the de Rham cohomology with coefficients in $M$ of a projective quotient $X_{\Gamma}$ of $X$ is that for the position of its Hodge filtration (e.g. due to $p$-adic Hodge theory its knowledge in case $d=1$ is another crucial point in \cite{br}); the second part of this paper is devoted to this question.
  
We discuss the content in more detail. Let now more generally $K$ be a non-Archimedean local field with ring of integers  ${\cal O}_K$ and residue field $k$. Let $M$ be a rational representation of $G={\rm GL}\sb {d+1}(K)$, i.e. a finite dimensional $K$-vector space $M$ together with a morphism of $K$-group varieties ${\rm GL}\sb {d+1}\to{\rm GL}(M)$. It is well known that for any compact open subgroup $H$ of $G$ there exists an $H$-stable free ${\mathcal O}_K$-module lattice in $M$; we fix one such choice $M^0$ for $H={\rm GL}\sb {d+1}({\mathcal O}_K)$. We choose a totally ramified extension $\dot{K}$ of $K$ of degree $d+1$ and twist the action of $G$ on $M\otimes_K\dot{K}$ by a suitable $\dot{K}$-valued character of $G$. We show that the choice of $M^0$ determines for any other maximal open compact subgroup $H\subset G$ a distinguished $H$-stable ${\mathcal O}_{\dot{K}}$-lattice in $M\otimes_K\dot{K}$ and the collection of these lattices can be assembled into a $G$-equivariant coefficient system on the Bruhat-Tits building ${\mathcal BT}$ of ${\rm PGL}\sb {d+1}$. In fact this is only a reinterpretation of our Proposition \ref{latteas}. We do not mention ${\mathcal BT}$ at all, we rather work with the $G$-equivariant semistable formal ${\mathcal O}_K$-scheme ${\mathfrak X}$ underlying Drinfel'd's symmetric space $X$ over $K$ of dimension $d+1$, as constructed in \cite{mus}. It is well known that the intersections of the irreducible components of ${\mathfrak X}\otimes k$ are in natural bijection with the simplices of ${\mathcal BT}$. Thus what we do is to construct from $M^0$ a constructible $G$-equivariant subsheaf ${\mathcal M}^0_{{\mathcal O}_{\dot{K}}}$ of the constant sheaf with value $M\otimes_K\dot{K}$ on ${\mathfrak X}$ such that ${\mathcal M}^0_{{\mathcal O}_{\dot{K}}}(U)$ for quasicompact open $U\subset {\mathfrak X}$ is an ${\mathcal O}_{\dot{K}}$-lattice in $M\otimes_K\dot{K}$.

We then consider the coherent ${\mathcal O}_{{\mathfrak X}}\otimes_{{\mathcal O}_{{K}}}{\mathcal O}_{\dot{K}}$-module sheaf ${\mathcal M}^0_{{\mathcal O}_{\dot{K}}}\otimes_{{\mathcal O}_{{K}}}{\mathcal O}_{{\mathfrak X}}$ and compute explicitly its reduction $({\mathcal M}^0_{{\mathcal O}_{\dot{K}}}\otimes_{{\mathcal O}_{{K}}}{\mathcal O}_{{\mathfrak X}})\otimes_{{\mathcal O}_{\dot{K}}}k$. See Theorem \ref{strumn} for our result. Similarly, let $\Omega^{\bullet}_{{\mathfrak X}}$ be the logarithmic de Rham complex of ${\mathfrak X}$ and let ${\mathcal Log}^s({\mathcal M}^0_{{\mathcal O}_{\dot{K}}})$ be the the $\pi$-adic completion of the subsheaf of ${\mathcal M}^0_{{\mathcal O}_{\dot{K}}}\otimes_{{\mathcal O}_K}\Omega^{s}_{{\mathfrak{X}}}$ consisting of logarithmic differential $s$-forms; we compute explicitly ${\mathcal Log}^s({\mathcal M}^0_{{\mathcal O}_{\dot{K}}})\otimes_{{\mathcal O}_{\dot{K}}}k$. See Theorem \ref{logred} for our result. 

As an application, assume now that $M|_{{\rm SL}\sb {d+1}(K)}$ is the trivial representation $K$, the standard representation $M=K^{d+1}$ or its dual $(K^{d+1})^*$. We show (Proposition \ref{globlo})\begin{gather}H^j({\mathfrak{X}},{\mathcal Log}^s({\mathcal M}^0_{{\mathcal O}_{\dot{K}}}))\cong H^j({\mathfrak{X}},{\mathcal M}^0_{{\mathcal O}_{\dot{K}}}\otimes_{{\mathcal O}_K}\Omega^{s}_{{\mathfrak{X}}})\label{introiso}\end{gather}for any $j$ and any $s$. Using the above computations the proof of (\ref{introiso}) is reduced to the statement that for any irreducible component $Y$ of ${\mathfrak X}\otimes k$ --- such a $Y$ is the successive blowing up of ${\mathbb P}_k^{d}$ in all $k$-linear subspaces --- with logarithmic de Rham complex $\Omega_Y^{\bullet}$ we have $H^j(Y,\Omega_Y^s)= 0$ if $j\ne0$, and $H^0(Y,\Omega_Y^s)$ consists of global logarithmic differential $s$-forms on $Y$. 

In the second part of this paper (sections \ref{norehose} and \ref{rehose}) we develop general criteria for conjectures of Schneider raised in \cite{schn}. Let $\Gamma\subset {\rm SL}\sb {d+1}(K)$ be a cocompact discrete (torsionfree) subgroup; thus the quotient $X_{\Gamma}=\Gamma\backslash X$ of $X$ is a projective $K$-scheme \cite{mus}. Let $M$ be a $K[\Gamma]$-module with $\dim_K(M)<\infty$. Using the $\Gamma$-action (induced from the $\Gamma$-action on $M$) on the constant local system  on $X$ generated by $M$ we get a local system ${\mathcal {M}}^{\Gamma}$ on ${X}_{\Gamma}$. The Hodge spectral sequence\begin{gather}E_1^{r,s}=H^s(X_{\Gamma},{\mathcal {M}}^{\Gamma}\otimes_K\Omega^{r}_{X_{\Gamma}} )\Rightarrow H^{r+s}(X_{\Gamma},{\mathcal {M}}^{\Gamma}\otimes_K\Omega^{\bullet}_{{X}_{\Gamma}}    )\label{exhodss}\end{gather}gives rise to the Hodge filtration$$H^d=H^{d}(X_{\Gamma},{\mathcal {M}}^{\Gamma}\otimes_K\Omega^{\bullet}_{{X}_{\Gamma}})=F^0_{H}\supset F^1_{H}\supset\ldots\supset F^{d+1}_{H}=0.$$If $\kara(K)=0$ Schneider \cite{schn} conjectures that a splitting of $F_H^{\bullet}$ is given by the covering filtration $F_{\Gamma}^{\bullet}$ of $H^d$ arising from the expression of $H^d$ through the $\Gamma$-group cohomology of $M\otimes_KH_{dR}^*(X)$. Concretely, he expects $H^d=F_H^{i+1}\oplus F_{\Gamma}^{d-i}$ for $0\le i\le d-1$.

If $M$ underlies a $K$-rational representation of $G$ (and $\kara(K)=0$), Schneider defines another sheaf complex, quasiisomorphic with ${\mathcal {M}}^{\Gamma}\otimes_K\Omega^{\bullet}_{{X}_{\Gamma}}$, hence again a corresponding Hodge filtration $F^{\bullet}_{red}$ on $H^d$. He then conjectures $F^{\bullet}_{red}=F_H^{\bullet}$, thus in particular he conjectures $H^d=F_{red}^{i+1}\oplus F_{\Gamma}^{d-i}$. The particular interest in this last decomposition is that combined with yet another conjecture from \cite{schn} --- the degeneration of the 'reduced'\, Hodge spectral sequence --- it would allow the computation of $\Gamma$-group cohomology spaces $H^*({\Gamma},D)$ for certain 'holomorphic discrete series representations'\, $D$ of $G$. 

For the trivial representation $M=K$ the conjectures were proven first by Iovita and Spiess \cite{iovspi}, later proofs were given by de Alon and Shalit (\cite{alde}, using harmonic analysis) and the author (\cite{hk}, using $p$-adic Hodge theory). The main tool in the approach of Iovita and Spiess is a certain subcomplex of $\Omega_X^{\bullet}(X)$ consisting of bounded logarithmic differential forms on $X$. For more general $M$ this complex does not seem to generalize well, essentially because there is no integral structure in the complex $M\otimes_K\Omega_X^{\bullet}(X)$ of {\it global} forms. This led us to consider a $K$-vector space sub{\it sheaf} complex $\mathbb{L}^{\bullet}(M)$ of $\underline{M}\otimes_K\Omega^{\bullet}_{X}$ on $X$ which should replace the global logarithmic differential forms. We show that the filtration $F_{\Gamma}^{\bullet}$ can be redefined in terms of $\mathbb{L}^{\bullet}(M)$ and obtain criteria for the above splitting conjectures and the degeneration of (\ref{exhodss}) which avoid $\Gamma$-group cohomology of global objects. A certain variant of $\mathbb{L}^{\bullet}(M)$, the  $K$-vector space sheaf complex $\mathbb{L}_D^{\bullet}(M)$, leads to a similar criterion for the splitting $H^d=F_{red}^{i+1}\oplus F_{\Gamma}^{d-i}$ and the degeneration of the 'reduced'\, Hodge spectral sequence. The general hope is that, working as indicated with integral (or bounded) structures inside $\mathbb{L}^{\bullet}(M)$ or $\mathbb{L}_D^{\bullet}(M)$, we can reduce to problems in characteristic $p$ and work {\it locally} on the reduction of the natural formal scheme underlying $X$. This approach worked out in \cite{mathan} in dimension $d=1$ where we used integral structures inside $\mathbb{L}_D^{\bullet}(M)$ to prove $H^1=F_{red}^{1}\oplus F_{\Gamma}^{1}$ and the degeneration conjecture. Here, as suggested above, we use integral structures inside $\mathbb{L}^{\bullet}(M)$ provided by the first part of this paper to prove (for arbitrary dimension $d$):\\

{\bf Theorem:} (see Corollary \ref{hosssp}, Theorem \ref{gammafi} and the remarks given there) {\it Suppose that $M|_{{\rm SL}\sb {d+1}(K)}=K$, $M|_{{\rm SL}\sb {d+1}(K)}=K^{d+1}$ or $M|_{{\rm SL}\sb {d+1}(K)}=(K^{d+1})^*$.\\(a) For arbitrary $\kara(K)$ the Hodge spectral sequence (\ref{exhodss}) degenerates in $E_1$. The Hodge filtration $F_H^{\bullet}$ has a canonical splitting defined through logarithmic differential forms.\\(b) If $\kara(K)=0$ we have $F_H^{\bullet}=F_{red}^{\bullet}$ and the splitting in (a) is given by the filtration $F_{\Gamma}^{\bullet}$:$$H^{d}(X_{\Gamma},{\mathcal {M}}^{\Gamma}\otimes_K\Omega^{\bullet}_{{X}_{\Gamma}}    )=F_H^{i+1}\bigoplus F_{\Gamma}^{d-i}\quad\quad(0\le i\le d-1).$$}

It seems that even for $M=K$ the degeneration in (a) in case $\kara(K)>0$ was unknown before. \\

{\it Notations:} We fix $d\in \mathbb{N}$ and {\it enumerate the rows and columns of} ${\rm GL}\sb {d+1}$-{\it elements by} $0,\ldots,d$. We denote by ${U}$ the subgroup of ${\rm GL}\sb {d+1}$ consisting of unipotent upper triangular matrices,$${U}=\{(a_{ij})_{0\le i,j\le d}\in{\rm GL}\sb {d+1}\quad|\quad\,a_{ii}=1\,\mbox{for all}\,i,\quad \,a_{ij}=0\,\,\mbox{if}\,\,i>j\}.$$For $r\in\mathbb{R}$ define $\lfloor r\rfloor, \lceil r\rceil\in\mathbb{Z}$ by requiring $\lfloor r\rfloor\le r<\lfloor r\rfloor+1$ and $\lceil r\rceil-1<r\le\lceil r\rceil$. For a divisor $D$ on an integral scheme $X$ we denote by ${\mathcal L}_X(D)$ the associated line bundle on $X$; we will always consider it as a subsheaf of the constant sheaf with value the function field of $X$. 

$K$ denotes a non-Archimedean local field, ${\cal O}_K$ its ring of integers, $\pi\in{\cal O}_K$ a fixed prime element and $k$ the residue field with $q$ elements, $q\in p^{\mathbb{N}}$. Let $\omega:K^{\times}_a\to\mathbb{Q}$ be the extension of the discrete valuation $\omega:K^{\times}\to\mathbb{Z}$ normalized by $\omega(\pi)=1$. We fix a totally ramified extension $\dot{K}=K(\dot{\pi})$ of $K$ with ring of integers ${\mathcal O}_{\dot{K}}$ such that $\dot{\pi}\in {\mathcal O}_{\dot{K}}$ satisfies $\dot{\pi}^{d+1}=\pi$.

We write $G={\rm GL}\sb {d+1}(K)$. Let $T$ be the torus of diagonal matrices in $G$ and let $X_*(T)$, resp. $X^*(T)$, denote the group of algebraic cocharacters, resp. characters, of $T$. For $0\le i\le d$ define the obvious cocharacter $e_i:{\mathbb{G}}_m\to{\rm GL}\sb {d+1}$, i.e. the one which sends $t$ to the diagonal matrix $(e_i(t))_{ij}$ with $e_i(t)_{ii}=t$, $e_i(t)_{jj}=1$ for $i\ne j$ and $e_i(t)_{j_1j_2}=0$ for $j_1\ne j_2$. The $e_i$ form a $\mathbb{R}$-basis of $X_*(T)\otimes{\mathbb{R}}$. The pairing $X_*(T)\times X^*(T)\to\mathbb{Z}$ which sends $(x,\mu)$ to the integer $\mu(x)$ such that $\mu(x(y))=y^{\mu(x)}$ for any $y\in\mathbb{G}_m$ extends to a duality between the $\mathbb{R}$-vector spaces $X_*(T)\otimes{\mathbb{R}}$ and $X^*(T)\otimes\mathbb{R}$. Let $\epsilon_0,\ldots,\epsilon_{d}\in X^*(T)$ denote the basis dual to $e_0,\ldots,e_{d}$. Let$$\Phi=\{\epsilon_i-\epsilon_j;\,\,0\le i,j\le d\,\,\mbox{and}\,i\ne j\}\quad\subset\quad X^*(T).$$

\section{Differential forms on rational varieties in characteristic $p>0$}
\label{cohmodp}

The action of ${\rm GL}_{d+1}(k)={\rm GL}(k^{d+1})$ on $(k^{d+1})^*=\Hom_k(k^{d+1},k)$ defines an action of ${\rm GL}_{d+1}(k)$ on the affine $k$-scheme associated with $(k^{d+1})^*$, and this action passes to an action of ${\rm GL}_{d+1}(k)$ on the projective space $$Y_0={\mathbb P}((k^{d+1})^*)\cong {\mathbb P}_k^{d}.$$For $0\le j\le d-1$ let ${\mathcal V}_0^j$ be the set of all $k$-rational linear subvarieties $Z$ of $Y_0$ with $\dim(Z)=j$, and let ${\mathcal V}_0=\bigcup_{j=0}^{d-1}{\mathcal V}_0^j$. The sequence of projective $k$-varieties$$Y=Y_{d-1}{\longrightarrow}Y_{d-2}{\longrightarrow}\ldots{\longrightarrow}Y_0$$is defined inductively by letting $Y_{j+1}\to Y_j$ be the blowing up of $Y_j$ in the strict transforms (in $Y_j$) of all $Z\in {\mathcal V}_0^j$. The set$${\mathcal V}=\mbox{the set of all strict transforms in}\,\,Y\,\,\mbox{of elements of}\,\,{\mathcal V}_0$$is a set of divisors on $Y$. The action of ${\rm GL}_{d+1}(k)$ on $Y_0$ naturally lifts to an action of ${\rm GL}\sb {d+1}(k)$ on $Y$. Let $\Xi_0,\ldots,\Xi_{d}$ be the standard projective coordinate functions on $Y_0$ and hence on $Y$ corresponding to the canonical basis of $(k^{d+1})^*$; hence $Y_0=\proj(k[\Xi_i;\,0\le i\le d])$. Denote by $\Omega^{\bullet}_{Y}$ the de Rham complex on $Y$ with logarithmic poles along the normal crossings divisor $\sum_{V\in{\mathcal V}}V$ on $Y$. For $i,j\in\{0,\ldots,d\}$ and $g\in {\rm GL}_{d+1}(k)$ we call$$g \dlog(\frac{\Xi_i}{\Xi_j})$$a logarithmic differential $1$-form on $Y$. We call an exterior product of logarithmic differential $1$-forms on $Y$ a logarithmic differential form on $Y$.

\begin{pro}\label{logdif} For each $0\le s\le d$ we have $H^t(Y,\Omega_Y^s)=0$ for all $t>0$. The $k$-vector space $H^0(Y,\Omega_Y^s)$ is the one generated by all logarithmic differential forms.\end{pro}

{\sc Proof:} In \cite{holdis} we derive this from a general vanishing theorem for higher cohomology of a certain class of line bundles on $Y$. Note that a corresponding statement over a field $F$ of characteristic zero is shown in \cite{iovspi} section 3: the de Rham cohomology of the complement of a finite set of $F$-rational hyperplanes in ${\mathbb P}^d_F$ is generated by (global) logarithmic differential forms. And the analoguous statement for the Monsky-Washnitzer cohomology of $Y^0=Y-\cup_{V\in{\mathcal V}}V$ was shown in \cite{ds}.\hfill$\Box$\\

{\it Remark:} In \cite{holdis} we give a $k$-basis for $H^0(Y,\Omega_Y^s)$ consisting of logarithmic differential forms as follows. For a subset $\tau\subset\{1,\ldots,d\}$ let $$U(k)(\tau)=\{(a_{ij})_{0\le i,j\le d}\in{{{U}}}(k)\quad|\quad a_{ij}=0\,\mbox{if}\, j\notin\{i\}\cup\tau \}.$$For $0\le s\le d$ denote by ${\mathcal P_s}$ the set of subsets of $\{1,\ldots,d\}$ consisting of $s$ elements. The following set is a $k$-basis of $H^0(Y,\Omega_Y^s)$:
$$\{A.\bigwedge_{t\in\tau}\dlog(\frac{\Xi_t}{\Xi_0})\quad|\quad \tau\in {\mathcal P_s}, A\in U(k)(\tau)\}.$$\hfill$\Box$\\

Let $D$ be a divisor on $Y$ of the type$$D=\sum_{V\in{\mathcal V}}b_{V}V$$with certain $b_{V}\in\mathbb{Z}$. We view ${\mathcal{L}}_{Y}(D)$ as a subsheaf of the constant sheaf $\underline{k(Y)}$ with value the function field $k(Y)$ of $Y$; hence we view $\Omega^{\bullet}_{Y}\otimes_{{\mathcal{O}}_{Y}}{\mathcal{L}}_{Y}(D)$ as a subsheaf of the constant sheaf with value the de Rham complex of $k(Y)/k$. The differential on the latter provides us with a differential on $\Omega^{\bullet}_{Y}\otimes_{{\mathcal{O}}_{Y}}{\mathcal{L}}_{Y}(D)$.

Consider the open and ${\rm GL}_{d+1}(k)$-stable subscheme $$Y^0=Y-\cup_{V\in{\mathcal V}}V$$of $Y$; let us write $$\iota:Y^0\to Y$$ for the embedding and $\Omega^{\bullet}_{Y^0}=\Omega^{\bullet}_{Y}|_{Y^0}$. 

For $0\le s\le d$ let ${\mathbb L}^{s}_Y$ be the $k$-vector subspace of $\Omega_Y^s(Y^0)$ generated by all $s$-forms $\eta$ of the type\begin{gather}\eta=y_1^{m_1}\cdots y_s^{m_s}\dlog(y_1)\wedge\ldots\wedge\dlog(y_s)\label{frofred}\end{gather}with $m_j\in\mathbb{Z}$ and $y_1,\ldots,y_d\in{\mathcal O}_Y^{\times}(Y^0)$ such that $y_j=\theta_j/\theta_0$ for a suitable (adapted to $\eta$) isomorphism of $k$-varieties $Y_0\cong\proj(k[\theta_j]_{0\le j\le d})$. From Proposition \ref{logdif} it follows that $H^0(Y,\Omega_Y^s)$ is the $k$-vector subspace of ${\mathbb L}^{s}_Y$ generated by all $s$-forms $\eta$ of type (\ref{frofred}) with $m_j=0$ for all $1\le j\le s$.

Let $\underline{\mathbb L}_Y^{s}$, resp. $\underline{\mathbb L}_Y^{s,0}$, be the constant sheaf on $Y$ with value ${\mathbb L}_Y^{s}$, resp. with value $H^0(Y,\Omega_Y^s)$. For a divisor $D$ as above we define$${\mathbb L}^{s}(D)=\underline{\mathbb L}_Y^{s}\quad\bigcap\quad {\mathcal L}_Y(D)\otimes\Omega_Y^{s},$$$${\mathbb L}^{s,0}(D)=\underline{\mathbb L}^{s,0}_{Y}\quad\bigcap\quad{\mathcal L}_{Y}(D)\otimes\Omega_Y^{s},$$the intersections taking place inside the push forward $\iota_*\Omega_{Y^0}^{s}$.

\begin{satz}\label{konskomop} (a) Suppose $b_V\in\{-1,0\}$ for all $V$. Then the inclusions ${\mathbb L}^{s,0}(D)\hookrightarrow{\mathbb L}^{s}(D)\hookrightarrow{\mathcal L}_Y(D)\otimes\Omega_Y^{s}$ induce for all $j$ isomorphisms$$H^j(Y,{\mathbb L}^{s,0}(D))\cong H^j(Y,{\mathbb L}^{s}(D))\cong H^j(Y,{\mathcal L}_{Y}(D)\otimes\Omega_Y^{s}).$$(b) Let $S$ be a non-empty subset of ${\mathcal V}$ such that $E=\cap_{V\in{S}}V$ is non-empty. Define the subsheaf ${\mathbb L}_E^s(0)$ of $\Omega^s_Y\otimes_{{\mathcal O}_Y}{\mathcal O}_E$ as the image of the composite ${\mathbb L}^s(0)\to\Omega^s_Y\to\Omega^s_Y\otimes_{{\mathcal O}_Y}{\mathcal O}_E$. Then the inclusion induces for all $j$ an isomorphism$$H^j(Y,{\mathbb L}_E^{s}(0))\cong H^j(Y,\Omega_Y^{s}\otimes_{{\mathcal O}_Y}{\mathcal O}_E).$$
\end{satz}

{\sc Proof:} (a) First we consider the case $D=0$, i.e. $b_V=0$ for all $V$. The sheaf ${\mathbb L}^{s,0}(0)$ is constant with value $H^0(Y,\Omega_Y^{s})$, hence we get $H^j(Y,{\mathbb L}^{s,0}(0))=H^j(Y,\Omega_Y^{s})$ for all $j$ from Proposition \ref{logdif}. In order to also compare with $H^j(Y,{\mathbb L}^{s}(0))$ choose a sequence $(\eta_n)_{n\ge1}$ of elements of ${\mathbb L}^{s}_Y$ of the form (\ref{frofred}) such that $\{\eta_n;\,n\ge1\}$ is a $k$-basis of ${\mathbb L}^{s}_Y/H^0(Y,\Omega_Y^s)$. For $n\ge0$ let $\underline{\mathbb L}_Y^{s,n}$ be the constant subsheaf of $\underline{\mathbb L}_Y^{s}$ on $Y$ generated over $k$ by $H^0(Y,\Omega_Y^s)$ and $\{\eta_{i};\,n\ge i\ge1\}$. Letting$${\mathbb L}^{s,n}(0)=\underline{\mathbb L}_Y^{s,n}\quad\bigcap \quad\Omega_Y^{s}$$we have $${\mathbb L}^{s}(0)=\bigcup_{n\ge0}{\mathbb L}^{s,n}(0)$$and since $Y$ is quasicompact (so that taking cohomology commutes with direct limits) it suffices to show$$H^j(Y,{\mathbb L}^{s,n}(0))=\left\{\begin{array}{l@{\quad:\quad}l}H^0(Y,\Omega_Y^s)&j=0\\0&j>0\end{array}\right.$$for all $n\ge0$. For $n=0$ we already did this, for $n>0$ it suffices, by induction, to show$$H^j(Y,\frac{{\mathbb L}^{s,n}(0)}{{\mathbb L}^{s,n-1}(0)})=0$$for all $j$. Let $W\subset Y$ be the maximal open subscheme on which the class of $\eta_n$ as a section of ${\mathbb L}^{s,n}(0)/{\mathbb L}^{s,n-1}(0)$ is defined. Thus if $\xi:W\to Y$ denotes the open embeddding and $\underline{k}_Y$ the constant sheaf on $Y$ with value $k$ then sending $1\in k$ to $\eta_n$ defines an isomorphism$$\xi_!\xi^{-1}\underline{k}_Y\cong\frac{{\mathbb L}^{s,n}(0)}{{\mathbb L}^{s,n-1}(0)}.$$If we had $W=Y$ then the induction hypothesis $H^1(Y,{\mathbb L}^{s,n-1}(0))=0$ and the long exact cohomology sequence associated with $$0\longrightarrow {\mathbb L}^{s,n-1}(0)\longrightarrow{\mathbb L}^{s,n}(0)\longrightarrow \frac{{\mathbb L}^{s,n}(0)}{{\mathbb L}^{s,n-1}(0)}\longrightarrow0$$would imply that there existed $a_1,\ldots,a_{n-1}\in k$ such that $\eta_n+\sum_{i=1}^{n-1}a_i\eta_i$ is a global section of ${\mathbb L}^{s,n}(0)$, in particular of $\Omega_Y^{s}$. But this would contradict the fact that $\{\eta_n;\,n\ge1\}$ is a $k$-basis of ${\mathbb L}^{s}_Y/H^0(Y,\Omega_Y^s)$. Hence $W\ne Y$. On the other hand we may write $\eta_n=y_1^{m_1}\cdots y_s^{m_s}\dlog(y_1)\wedge\ldots\wedge\dlog(y_s)$ with $y_j=\theta_j/\theta_0$ as in (\ref{frofred}) and it is clear that $C=Y-W$ is the pull back under $Y\to Y_0$ of a union of some hyperplanes $\mathbb{V}(\theta_i)\subset Y_0$. In particular $C$ is connected. Denote by $\gamma:C\to Y$ the closed embedding. The long exact cohomology sequence associated with$$0\longrightarrow\xi_!\xi^{-1}\underline{k}_Y\longrightarrow\underline{k}_{Y}\longrightarrow\gamma_*\gamma^{-1}\underline{k}_{Y}\longrightarrow0$$shows $H^j(Y,\xi_!\xi^{-1}\underline{k}_{Y})=0$ for all $j$ because $C$ is non empty and connected. The induction and thus the discussion of the case $D=0$ is finished.
 
To treat arbitrary $D$ with $b_V\in\{-1,0\}$ for all $V$ we induct on $\dim(Y)$ and on $r(D)=\sum_{V\in{\mathcal{V}}}|b_V|$. We will only show $H^j(Y,{\mathbb L}^{s}(D))= H^j(Y,{\mathcal L}_{Y}(D)\otimes_{{\mathcal O}_Y}\Omega_Y^{s})$ (which is the relevant statement for the subsequent sections), the proof of $H^j(Y,{\mathbb L}^{s,0}(D))= H^j(Y,{\mathcal L}_{Y}(D)\otimes_{{\mathcal O}_Y}\Omega_Y^{s})$ is literally the same (replace each occurence of ${\mathbb L}^{s}(D)$ with ${\mathbb L}^{s,0}(D)$).

Assume $b_V=-1$ for some $V$. Let $D'=D+V$. We want to compare the exact sequences
$$0\longrightarrow{\mathbb L}^{s}(D)\longrightarrow{\mathbb L}^{s}(D')\longrightarrow{\mathbb L}^{s}_{V}(D')\longrightarrow0$$(the sheaf ${\mathbb L}^{s}_{V}(D')$ being defined such that this is an exact sequence) and$$0\longrightarrow{\mathcal L}_Y(D)\otimes\Omega_Y^{s}\longrightarrow{\mathcal L}_Y(D')\otimes\Omega_Y^{s}\longrightarrow{\mathcal L}_Y(D')\otimes\Omega_Y^{s}\otimes{\mathcal O}_V\longrightarrow0.$$Since $r(D')<r(D)$ the induction hypothesis says that the map between the respective second terms induces isomorphisms in cohomology. It will be enough to prove the same for the respective third terms. From \cite{ito} (see also \cite{holdis}) it follows that $V$ decomposes as$$V=Y^1\times Y^2$$such that both $Y^t$ are successive blowing ups of projective spaces of dimensions smaller than $d$ in all $k$-rational linear subvarieties, just as $Y$ is the successive blowing up of projective space of dimension $d$ in all $k$-rational linear subvarieties. Denote by ${\mathcal V}^t$ the corresponding set of divisors on $Y^t$ (like the set ${\mathcal V}$ of divisors on $Y$) and let $\Omega^{\bullet}_{Y^t}$ denote the logarithmic de Rham complex on $Y^t$ with logarithmic poles along ${\mathcal V}^t$. Let $\Omega^{\bullet}_V$ denote the logarithmic de Rham complex on $V$ with logarithmic poles along all divisors which are pull backs of elements of ${\mathcal V}^1$ or ${\mathcal V}^2$. Then$$\Omega^{\bullet}_V=\Omega^{\bullet}_{Y^1}\otimes_k\Omega^{\bullet}_{Y^2}.$$Let $D_V$ be the divisor on $V$ induced by $D$. More precisely, $D_V=\sum b_{W}(W\cap V)$, the sum ranging over all $W\in\mathcal{V}$ which intersect $V$ transversally. It also follows from \cite{ito} (and \cite{holdis}) that $D_V$ is of the form $D_{Y^1}^V+D_{Y^2}^V$ where $D_{Y^t}^V$ for $t=1,2$ is the pull back to $V$, via the projection $V\to Y^t$, of a divisor $D_{Y^t}$ on $Y^t$ which is the sum, with multiplicities in $\{0,-1\}$, of elements of ${\mathcal V}^t$. The above then generalizes as\begin{gather}{\mathcal L}_V(D_V)\otimes\Omega^{\bullet}_V=({\mathcal L}_{Y^1}(D_{Y^1})\otimes\Omega^{\bullet}_{Y^1})\otimes_k{(\mathcal L}_{Y^2}(D_{Y^2})\otimes\Omega^{\bullet}_{Y^2}).\label{lokue}\end{gather}Choose an isomorphism $Y_0\cong\proj(k[\theta_j]_{0\le j\le d})$ and elements $0\le j_1\ne j_2\le d$ such that $y=\theta_{j_1}/\theta_{j_2}\in{\mathcal O}_Y(U)$ is an equation for $V\cap U$ in a suitable open subset $U$ of $Y$ with $V\cap U\ne\emptyset$. We have an exact sequence\begin{gather}0\longrightarrow{\mathcal L}_V(D_V)\otimes\Omega^{s-1}_V\stackrel{\wedge\dlog(y)}{\longrightarrow}{\mathcal L}_Y(D')\otimes\Omega_Y^{s}\otimes{\mathcal O}_V\longrightarrow{\mathcal L}_V(D_V)\otimes\Omega^{s}_V\longrightarrow0\label{fullex}\end{gather}where by (\ref{lokue}) the extreme terms (take $s'=s$ and $s'=s-1$) decompose as 
\begin{gather}{\mathcal L}_V(D_V)\otimes\Omega^{s'}_V\cong\bigoplus_{s_1+s_2=s'}({\mathcal L}_{Y^1}(D_{Y^1})\otimes\Omega^{s_1}_{Y^1})\otimes_k({\mathcal L}_{Y^2}(D_{Y^2})\otimes\Omega^{s_2}_{Y^2}).\label{dide}\end{gather}On the other hand, define for $t=1,2$ the sheaves ${\mathbb L}^{\bullet}(D_{Y^t})$ on $Y^t$ just as we defined the sheaves ${\mathbb L}^{\bullet}(.)$ on $Y$ (and use the same name for their push forward to $Y$). Then using the decomposition (\ref{dide}) we may view the sheaf$$\bigoplus_{s_1+s_2=s'}{\mathbb L}^{s_1}(D_{Y^1})\otimes_k{\mathbb L}^{s_2}(D_{Y^2})$$as a subsheaf of ${\mathcal L}_V(D_V)\otimes\Omega^{s'}_V$ and a local consideration shows that (\ref{fullex}) restricts to an exact sequence
\begin{gather}0\longrightarrow\bigoplus_{s_1+s_2=s-1}{\mathbb L}^{s_1}(D_{Y^1})\otimes_k{\mathbb L}^{s_2}(D_{Y^2})\stackrel{\wedge\dlog(y)}{\longrightarrow}{\mathbb L}^{s}_{V}(D')\longrightarrow\bigoplus_{s_1+s_2=s}{\mathbb L}^{s_1}(D_{Y^1})\otimes_k{\mathbb L}^{s_2}(D_{Y^2}){\longrightarrow}0.\label{lex}\end{gather}Comparing the long exact cohomology sequences associated with (\ref{fullex}) and (\ref{lex}) we conclude that to show that$$H^j(Y,{\mathbb L}^{s}_{V}(D'))\longrightarrow H^j(Y,{\mathcal L}_Y(D')\otimes\Omega_Y^{s}\otimes{\mathcal O}_V)$$is an isomorphism for any $j$, it suffices to show that$$H^j(Y,\bigoplus_{s_1+s_2=s'}{\mathbb L}^{s_1}(D_{Y^1})\otimes_k{\mathbb L}^{s_2}(D_{Y^2}))\longrightarrow H^j(Y,\bigoplus_{s_1+s_2=s'}({\mathcal L}_{Y^1}(D_{Y^1})\otimes\Omega^{s_1}_{Y^1})\otimes_k({\mathcal L}_{Y^2}(D_{Y^2})\otimes\Omega^{s_2}_{Y^2}))$$is an isomorphism, for $s'=s$ and $s'=s-1$. By the K\"unneth formula this reduces to showing that$$H^j(Y^t,{\mathbb L}^{s''}(D_{Y^t}))\longrightarrow H^j(Y^t,({\mathcal L}_{Y^t}(D_{Y^t})\otimes\Omega^{s''}_{Y^t}))$$is an isomorphism, for any $s''$ and $t\in\{1,2\}$. But this follows from our induction hypothesis since the dimension of $Y^t$ is smaller than that of $Y$.\\(b) We have an exact sequence$$0\longrightarrow{\mathcal L}_Y(-\sum_{V\in S}V)\longrightarrow\bigoplus_{T\subset S\atop|T|=|S|-1}{\mathcal L}_Y(-\sum_{V\in T}V)\longrightarrow\ldots\longrightarrow\bigoplus_{V\in S}{\mathcal L}_Y(-V)\longrightarrow{\mathcal O}_Y\longrightarrow{\mathcal O}_E\longrightarrow0$$which yields a similar exact sequence by tensoring with $\Omega_Y^s$. A local consideration shows that the latter exact sequence restricts to an exact sequence$$0\longrightarrow{\mathbb L}^s(-\sum_{V\in S}V)\longrightarrow\bigoplus_{T\subset S\atop|T|=|S|-1}{\mathbb L}^s(-\sum_{V\in T}V)\longrightarrow\ldots\longrightarrow\bigoplus_{V\in S}{\mathbb L}^s(-V)\longrightarrow{\mathbb L}^s(0)\longrightarrow{\mathbb L}_E^s(0)\longrightarrow0.$$It follows that it is enough to show that for all subsets $T\subset S$ and any $j$ the map$$H^j(Y,{\mathbb L}^s(-\sum_{V\in T}V))\longrightarrow H^j(Y,{\mathcal L}_Y(-\sum_{V\in T}V)\otimes_{{\mathcal O}_Y}\Omega_Y^{s})$$is an isomorphism. But this follows from part (a).\hfill$\Box$\\

{\it Remark:} (not needed in the sequel) If $-1\le b_V\le p-1$ for all $V$ then the inclusion ${\mathbb L}^{\bullet,0}(D)\hookrightarrow{\mathcal L}_Y(D)\otimes\Omega_Y^{\bullet}$ induces isomorphisms$$H^j(Y,{\mathbb L}^{\bullet,0}(D))\cong H^j(Y,{\mathcal L}_Y(D)\otimes\Omega_Y^{\bullet}).$$To see this let $D'=\sum_{V\in{\mathcal V}}b'_{V}V$ with $b'_V=\min\{0,b_V\}$. Then Theorem \ref{konskomop} applies to $D'$. Now note that on the one hand ${\mathbb L}^{\bullet,0}(D')={\mathbb L}^{\bullet,0}(D)$ (logarithmic differential forms have pole orders at most one) and on the other hand ${\mathcal L}_Y(D')\otimes\Omega_Y^{\bullet}\longrightarrow{\mathcal L}_Y(D)\otimes\Omega_Y^{\bullet}$ is a quasiismorphism (use that any $b_V>0$ is invertible in $k$).

\section{Reduction of rational $G$-representations}
\label{rera}

Let $\overline{T}=T/K^{\times}$. For $\mu=\sum_{j=0}^da_j\epsilon_j\in X^*(T)$ let \begin{gather}\overline{\mu}=(\frac{1}{d+1}\sum_{j=0}^da_j)(\sum_{j=0}^d\epsilon_j)-\mu,\label{mubar}\end{gather}an element of the subspace $X^*(\overline{T})\otimes\frac{1}{d+1}.\mathbb{Z}$ of $X^*({T})\otimes\frac{1}{d+1}.\mathbb{Z}$. If for $0\le j\le d$ we let\begin{gather}\overline{a}_{j}(\mu)=\frac{(\sum_{i\ne j}a_i)-da_{j}}{d+1}\label{aquer}\end{gather}then$$\overline{\mu}=\sum_{j=0}^d\overline{a}_j(\mu)\epsilon_j.$$

Let $M$ be an irreducible $K$-rational representation of $G$. For a weight $\mu\in X^{*}(T)$ let $M_{\mu}$ be the maximal subspace of $M$ on which $T$ acts through $\mu$. 

\begin{lem}\label{dichot} The number $$|M|=\sum_{i=0}^da_{i}$$ for $\mu=\sum_{i=0}^da_{i}\epsilon_i\in X^*(T)$ such that $M_{\mu}\ne 0$ is independent of the choice of such a $\mu$; for such $\mu$ we have $\overline{\mu}\in X^*(\overline{T})$ if and only if $|M|\in(d+1).\mathbb{Z}$, if and only if there is a $h\in\mathbb{Z}$ such that the center of $G$ acts trivially on $M\otimes_K\det^h$. 
\end{lem}

{\sc Proof:} This is clear since all $\mu$ with $M_{\mu}\ne 0$ differ by linear combinations of elements of $\Phi$ (see \cite{jan} II.2.2).\hfill$\Box$\\

We fix a ${\rm GL}\sb {d+1}/{\mathcal O}_K$-invariant ${\mathcal O}_K$-lattice $M^0$ in $M$ (see \cite{jan} I.10.4). 

\begin{lem}\label{intspan} We have $M^0=\bigoplus_{\mu\in X^*(T)}M^0_{\mu}$ with $M_{\mu}^0=M^0\cap M_{\mu}$.
\end{lem}

{\sc Proof:} We reproduce a proof from notes of Schneider and Teitelbaum. Fix $\mu\in X^*(T)$. It suffices to construct an element $\Pi_{\mu}$ in the algebra of distributions $\mbox{Dist}({\rm GL}\sb {d+1}/\mathbb{Z})$ (i.e. defined over $\mathbb{Z}$) which on $M$ acts as a projector onto $M_{\mu}$. For $0\le i\le d$ let $H_i=(de_i)(1)\in\mbox{Lie}({\rm GL}\sb {d+1}/\mathbb{Z})$; then $d\mu'(H_i)\in\mathbb{Z}$ (inside $\mbox{Lie}(\mathbb{G}_m/\mathbb{Z})$) for any $\mu'\in X^*(T)$. According to \cite{hum} Lemma 27.1 we therefore find a polynomial $\Pi\in{\mathbb{Q}}[X_0,\ldots,X_d]$ such that $\Pi(\mathbb{Z}^{d+1})\subset \mathbb{Z}$, $\Pi(d\mu(H_0),\ldots,d\mu(H_d))=1$ and $\Pi(d\mu'(H_0),\ldots,d\mu'(H_d))=0$ for any $\mu'\in X^*(T)$ such that $\mu'\ne\mu$ and $M_{\mu'}\ne0$. Moreover \cite{hum} Lemma 26.1 says that $\Pi$ is a $\mathbb{Z}$-linear combination of polynomials of the form $${X_0\choose b_0}\cdots{X_d\choose b_d}\quad\mbox{with integers}\quad b_0,\ldots,b_d\ge0.$$Thus \cite{jan} II.1.12 implies that $$\Pi_{\mu}=\Pi(H_0,\ldots,H_d)$$lies in $\mbox{Dist}({\rm GL}\sb {d+1}/\mathbb{Z})$. By construction it acts on $M$ as a projector onto $M_{\mu}$.\hfill$\Box$\\

We return to the setting from section \ref{cohmodp}. For $\emptyset\ne\sigma\subsetneq\{0,\ldots,d\}$ denote by $V^0_{\sigma}$ the common zero set in $Y_0$ of all $\Xi_j$ with $j\in\sigma$, and let $V_{\sigma}\in{\mathcal V}$ be its strict transform under $Y\to Y_0$. Denote by $Y'$ the open subscheme of $Y$ obtained by deleting all divisors $V\in{\mathcal V}$ which are {\it not} of the particular form $V=V_{\sigma}$ for some $\emptyset\ne\sigma\subsetneq\{0,\ldots,d\}$. Then $Y^0\subset Y'\subset Y$ and $U(k)$ acts on $Y$ and on $Y^0$ and moreover $U(k).Y'=Y$ (the $U(k)$-translates of $Y'$ cover $Y$). For each $V\in{\mathcal V}$ there is a unique $\emptyset\ne\sigma\subsetneq\{0,\ldots,d\}$ such that there exists a $g\in U(k)$ with $gV_{\sigma}=V$, see \cite{holdis}.

Let $\widetilde{M}=M^0/\pi.M^0$. The decomposition from Lemma \ref{intspan} induces a corresponding decomposition $$\widetilde{M}=\bigoplus_{\mu\in X^*(T)}\widetilde{M}_{\mu}.$$Denote again by ${\widetilde{M}}$ the constant sheaf on $Y$ with value $\widetilde{M}$. We define a subsheaf $\widetilde{\mathcal M}[{\mathcal O}_{Y'}]$ of ${\widetilde{M}}\otimes_k\iota_*{\mathcal O}_{Y^0}|_{Y'}$ on $Y'$ by\begin{gather}\widetilde{\mathcal M}[{\mathcal O}_{Y'}]=\bigoplus_{\mu\in X^*(T)}\widetilde{M}_{\mu}\otimes_k{\mathcal L}_{Y'}(\sum_{\emptyset\ne\sigma\subsetneq\{0,\ldots,d\}}-\lceil\sum_{j\in\sigma}\overline{a}_j(\mu)\rceil (V_{\sigma}\cap Y')).\label{ystrdf}\end{gather}

\begin{lem}\label{redspread} $\widetilde{\mathcal M}[{\mathcal O}_{Y'}]$ extends uniquely to a ${\rm GL}_{d+1}(k)$-stable subsheaf $\widetilde{\mathcal M}[{\mathcal O}_{Y}]$ of ${\widetilde{M}}\otimes_k\iota_*{\mathcal O}_{Y^0}$. 
\end{lem}

{\sc Proof:} This can be checked directly, an easier variant of the proof of Theorem \ref{strumn} below. However, it is even a {\it consequence} of Theorem \ref{strumn}: explicitly,$$\widetilde{\mathcal M}[{\mathcal O}_{Y}]=\frac{({\mathcal M}^0_{{\mathcal O}_{\dot{K}}}\otimes_{{\mathcal O}_K}{\mathcal O}_{\mathfrak{X}})\otimes_{{\mathcal O}_{\dot{\mathfrak{X}}}}{\mathcal O}_Y}{{\mathcal O}_Y\mbox{-torsion}}$$in the notations used there.\hfill$\Box$\\

{\bf Definition:} We say that the weights of $M$ are small if for any $\mu\in X^*(T)$ with $M_{\mu}\ne0$ and for any $\emptyset\ne \sigma\subsetneq\{0,\ldots,d\}$ we have\begin{gather}0\le\lceil\sum_{j\in\sigma}\overline{a}_j(\mu)\rceil\le1.\label{smaldef}\end{gather}

\begin{lem}\label{wesm} The weights of $M$ are small if and only if, when regarded as a representation of ${\rm SL}_{d+1}(K)$, it is one of the following: the trivial representation, the standard representation $K^{d+1}$, or the dual $(K^{d+1})^*$ of the standard representation of ${\rm SL}_{d+1}(K)$.\hfill$\Box$\\
\end{lem} 

{\sc Proof:} One easily checks that $\mu=\sum_{i=0}^da_{i}\epsilon_i\in X^*(T)$ satisfies inequality (\ref{smaldef}) for all $\sigma$ if and only if all coefficients $a_i$ are the same (case (i)) or if there is precisely one $0\le i\le d$ with $a_i=a_j+1$ for all $j\ne i$ (case (ii)) or with $a_i=a_j-1$ for all $j\ne i$ (case (iii)). If $M|_{{\rm SL}_{d+1}(K)}=K$ the only weight occuring is as in case (i), if $M|_{{\rm SL}_{d+1}(K)}=K^{d+1}$ the weights occuring are as in case (ii), if $M|_{{\rm SL}_{d+1}(K)}=(K^{d+1})^*$ the weights occuring are as in case (iii).\hfill$\Box$\\

For $0\le s\le d$ consider the following sheaf $\widetilde{\mathcal M}[{\mathbb L}^{s}_Y]$ on $Y$:$$\widetilde{\mathcal M}[{\mathbb L}^{s}_Y]={\widetilde{M}}\otimes_k\underline{\mathbb L}_Y^{s}\quad\bigcap\quad\widetilde{\mathcal M}[{\mathcal O}_Y]\otimes_{{\mathcal O}_Y}\Omega_Y^{s},$$the intersection taking place inside $\iota_*({\widetilde{M}}\otimes_k\Omega_Y^{s}|_{Y^0})$.

\begin{satz}\label{klewe} If the weights of $M$ are small then the inclusion $\widetilde{\mathcal M}[{\mathbb L}^{s}_Y]\to\widetilde{\mathcal M}[{\mathcal O}_Y]\otimes_{{\mathcal O}_Y}\Omega_Y^{s}$ of sheaves on $Y$ induces isomorphisms$$H^*(Y,\widetilde{\mathcal M}[{\mathbb L}^{s}_Y])\cong H^*(Y,\widetilde{\mathcal M}[{\mathcal O}_Y]\otimes_{{\mathcal O}_Y}\Omega_Y^{s}).$$
\end{satz}

{\sc Proof:} Consider the following ordering on $X^*(T)$: define\begin{gather}\sum_{i=0}^da_i\epsilon_i>\sum_{i=0}^da'_i\epsilon_i\label{ordwei}\end{gather}if and only if there exists a $0\le i_0\le d$ such that $a_i=a_i'$ for all $i<i_0$, and $a_{i_0}>a'_{i_0}$. By \cite{jan} II.1.19 the filtration $(F^{\mu}M)_{\mu\in X^*(T)}$ of $M$ defined by\begin{gather}F^{\mu}M=\sum_{\mu'\in X^*(T)\atop \mu'\ge\mu}M_{\mu'}\label{weifiab}\end{gather}is stable for the action of ${U}(K)$. Hence the filtration $(F^{\mu}M^0)_{\mu\in X^*(T)}$ of $M^0$ defined by$$F^{\mu}M^0=\sum_{\mu'\in X^*(T)\atop \mu'\ge\mu}M^0_{\mu'}=F^{\mu}M\bigcap M^0$$is stable for the action of ${U}({\mathcal O}_K)$, and the induced filtrations $(F^{\mu}\widetilde{\mathcal M}[{\mathcal O}_Y])_{\mu\in X^*(T)}$ of $\widetilde{\mathcal M}[{\mathcal O}_Y]$ and $(F^{\mu}\widetilde{\mathcal M}[{\mathbb L}^{s}_Y])_{\mu\in X^*(T)}$ of $\widetilde{\mathcal M}[{\mathbb L}^{s}_Y]$ are $U(k)$-stable. We denote by $$Gr^{\mu}(.)=\frac{F^{\mu}(.)}{\sum_{\mu'>\mu}F^{\mu'}(.)}$$the respective graded pieces. To prove Theorem \ref{klewe} it is enough to show that for all $\mu\in X^*(T)$ the maps$$H^*(Y,Gr^{\mu}\widetilde{\mathcal M}[{\mathbb L}^{s}_Y])\longrightarrow H^*(Y,Gr^{\mu}\widetilde{\mathcal M}[{\mathcal O}_Y]\otimes_{{\mathcal O}_Y}\Omega_Y^{s})$$are isomorphisms. By definition (\ref{ystrdf}), the restriction of $\widetilde{\mathcal M}[{\mathcal O}_Y]$ to $Y'$ comes with a canonical splitting of the filtration $(F^{\mu}\widetilde{\mathcal M}[{\mathcal O}_Y])_{\mu\in X^*(T)}$, and this splitting shows $$Gr^{\mu}\widetilde{\mathcal M}[{\mathcal O}_Y]\otimes_{{\mathcal O}_Y}\Omega_Y^{s}|_{Y'}\cong {\widetilde{M}}_{\mu}\otimes_k{\mathcal L}_{Y'}(\sum_{\emptyset\ne\sigma\subsetneq\{0,\ldots,d\}}-\lceil\sum_{j\in\sigma}\overline{a}_j(\mu)\rceil (V_{\sigma}\cap Y'))\otimes_{{\mathcal O}_Y}\Omega^s_Y|_{Y'}.$$Moreover, the subsheaf $Gr^{\mu}\widetilde{\mathcal M}[{\mathbb L}^{s}_Y]|_{Y'}$ of $Gr^{\mu}\widetilde{\mathcal M}[{\mathcal O}_Y]\otimes_{{\mathcal O}_Y}\Omega_Y^{s}|_{Y'}$ can be identified with $$\widetilde{M}_{\mu}\otimes_k{\mathbb L}^{s}(\sum_{\emptyset\ne\sigma\subsetneq\{0,\ldots,d\}}-\lceil\sum_{j\in\sigma}\overline{a}_j(\mu)\rceil V_{\sigma})|_{Y'}.$$Thus, by $U(k)$-equivariance and since $U(k)Y'=Y$, the inclusion $Gr^{\mu}\widetilde{\mathcal M}[{\mathbb L}^{s}_Y]\to Gr^{\mu}\widetilde{\mathcal M}[{\mathcal O}_Y]\otimes_{{\mathcal O}_Y}\Omega_Y^{s}$ is of the form considered in Theorem \ref{konskomop}, tensored with (the constant sheaf generated by) $\widetilde{M}_{\mu}$. Hence we may conclude by Theorem \ref{konskomop}.\hfill$\Box$

\section{Sheaves of integral structures in $K[G]$-modules}

\label{secints}

The action of $G={\rm GL}_{d+1}(K)={\rm GL}(K^{d+1})$ on $(K^{d+1})^*=\Hom_K(K^{d+1},K)$ defines an action of $G$ on the affine $K$-scheme associated with $(K^{d+1})^*$, and this action passes to an action of $G$ on the projective space ${\mathbb P}((K^{d+1})^*)$. Drinfel'd's symmetric space $X$ is the $K$-rigid space$$X={\mathbb P}((K^{d+1})^*)-(\mbox{the union of all $K$-rational hyperplanes}).$$Clearly $X$ is stable for the action of $G$. Let ${\mathfrak{X}}$ be the strictly semistable formal ${\mathcal O}_K$-scheme with generic fibre $X$ introduced in \cite{mus}. Instead of recalling its formal definition here we recall its basic properties relevant for us. ${\mathfrak{X}}$ is covered by Zariski open subschemes which admit open immersions into the $\pi$-adic formal completion of $\spec({\mathcal O}_K[X_0,\ldots,X_d]/(X_0\ldots X_d-\pi))$. Each irreducible component of the reduction ${\mathfrak X}\otimes k$ of ${\mathfrak X}$ is isomorphic to the $k$-scheme $Y$ studied in the previous sections. The set of all irreducible components of ${\mathfrak X}\otimes k$ is in natural bijection with the set of vertices of the Bruhat Tits building of ${\rm PGL}\sb {d+1}/K$. More generally, if for $j\ge 0$ we let $F^j$ denote the set of non-empty intersections of $(j+1)$-many pairwise distinct irreducible components of ${\mathfrak{X}}\otimes k$, then $F^j$ is in natural bijection with the set of $j$-simplices of the Bruhat Tits building of ${\rm PGL}\sb {d+1}/K$. This bijection is $G$-equivariant for the natural extension of the action of $G$ on $X$ to an action of $G$ on ${\mathfrak{X}}$. We denote by $Y$ the central irreducible component of ${\mathfrak X}\otimes k$, i.e. the irreducible component of ${\mathfrak X}\otimes k$ which is characterized by the fact that the subgroup $K^{\times}.{\rm GL}\sb {d+1}({{\mathcal O}_K})$ of $G$ is the stabilizer of $Y$ (for the action of $G$ on the set $F^0$ of irreducible components of ${\mathfrak X}\otimes k$). We identify this $k$-scheme $Y$ (with its ${\rm GL}\sb {d+1}({{\mathcal O}_K}$)-action) with the one from the previous sections. We define the subset $$F^0_{A}=T.Y$$of $F^0$, the orbit of $Y\in F^0$ for the action of $T$ on $F^0$. (This set corresponds to the set of vertices in the standard apartment of the Bruhat Tits building of ${\rm PGL}\sb {d+1}/K$.) For $Z\in F^0_{A}$ and $\mu\in X^*(T)$ we let $\overline{\mu}\in X^*(\overline{T})\otimes\frac{1}{d+1}{\mathbb{Z}}$ as before and define $\overline{\mu}(Z)\in\frac{1}{d+1}.\mathbb{Z}$ as$$\overline{\mu}(Z)=-\omega(\overline{\mu}(t))\quad\mbox{with}\,\,t\in T\,\,\mbox{such that}\,\,t.Y=Z.$$For $Z\in F^0$ let ${\mathfrak{U}}_Z$ be the maximal {\it open} formal subscheme of ${\mathfrak X}$ such that ${\mathfrak{U}}_Z\otimes_{{\mathcal O}_K}k$ is contained in $Z$. For example, the indicated identification of the central irreducible component $Y$ of ${\mathfrak X}\otimes k$ with the $k$-scheme $Y$ from the previous sections restricts to an identification of open subschemes ${\mathfrak{U}}_Y\otimes_{{\mathcal O}_K}k=Y^0$ (with $Y^0\subset Y$ as defined in the previous sections). Also note that the union $\cup_{Z\in F^0}{\mathfrak{U}}_Z$ is disjoint in ${\mathfrak X}$ and that the closed points of $\cup_{Z\in F^0}{\mathfrak{U}}_Z\otimes k$ are exactly the non-singular closed points of the $k$-scheme ${\mathfrak X}\otimes k$.\\   

Let again $M$ be an irreducible $K$-rational representation of ${\rm GL}\sb {d+1}$ and fix a ${\rm GL}\sb {d+1}/{\mathcal O}_K$-invariant ${\mathcal O}_K$-lattice $M^0$ in $M$. Define the character $$\chi:G\to \dot{K}^{\times},\quad\quad g\mapsto\dot{\pi}^{-\omega(\det(g))}$$and let $G$ act on ${M}\otimes_K{\dot{K}}$ by multiplying with $\chi^{|M|}$ the (scalar extension $K\to\dot{K}$ of the) already given action of $G$ on ${M}$. The point of this twisting is that the ${\mathcal O}_{\dot{K}}$-submodule $M^0\otimes_{{\mathcal O}_K}{\mathcal O}_{\dot{K}}$ of ${M}\otimes_K{\dot{K}}$ is now stable not just for ${\rm GL}\sb {d+1}({{\mathcal O}_K})$ but even for the full stabilizer $K^{\times}.{\rm GL}\sb {d+1}({{\mathcal O}_K})$ of $Y$ in ${\mathfrak X}$. Of course, if $|M|\in(d+1).\mathbb{Z}$ then we could replace the above twisting by a twisting with a suitable power of the determinant character of $G$, and the base extension $K\to\dot{K}$ here and in the whole construction below could be avoided.

Let $\underline{M}_{\dot{K}}$ be the constant sheaf on ${\mathfrak{X}}$ with value $M\otimes_K\dot{K}$. The above action of $G$ on $M\otimes_K\dot{K}$ makes $\underline{M}_{\dot{K}}$ into a $G$-equivariant sheaf on ${\mathfrak{X}}$. Define $M_{\mu}^0$ as in Lemma \ref{intspan}. For $\mu\in X^{*}(T)$ let $\underline{M}^0_{\mu,{\mathcal{O}}_{\dot{K}}}$ be the constant subsheaf of $\underline{M}_{\dot{K}}$ with value $M^0_{\mu}\otimes_{{\mathcal O}_K}{\mathcal O}_{\dot{K}}$. For $Z\in F^0_{A}$ let \begin{gather}{\mathcal M}^0_{{\mathcal O}_{\dot{K}}}|_{\mathfrak{U}_{Z}}=\bigoplus_{\mu\in X^{*}(T)}\dot{\pi}^{(d+1)\overline{\mu}(Z)}\underline{M}^0_{\mu,{\mathcal{O}}_{\dot{K}}}|_{\mathfrak{U}_{Z}}.\label{numdef}\end{gather}

\begin{pro}\label{latteas} Formula (\ref{numdef}) (for all $Z\in F^0_{A}$) defines a subsheaf $${\mathcal M}^0_{{\mathcal O}_{\dot{K}}}|_{\bigcup_{Z\in F_A^0}\mathfrak{U}_{Z}}\subset\underline{M}_{\dot{K}}|_{\bigcup_{Z\in F_A^0}\mathfrak{U}_Z}.$$It extends to a $G$-stable subsheaf ${\mathcal M}^0_{{\mathcal O}_{\dot{K}}}$ of $\underline{M}_{\dot{K}}$ in finitely generated ${\mathcal O}_{\dot{K}}$-modules such that $${\mathcal M}^0_{{\mathcal O}_{\dot{K}}}\otimes_{{\mathcal O}_{\dot{K}}}\dot{K}=\underline{M}_{\dot{K}}.$$
\end{pro}

{\sc Proof:} (Here we benefited from notes of Schneider and Teitelbaum.) We need some more notations. For $0\le i,j\le d$ and $i\ne j$ consider the morphism of algebraic groups over $\mathbb{Z}$\begin{gather}\widetilde{\alpha}_{ij}:\mathbb{G}_a\longrightarrow {\rm GL}\sb {d+1},\quad u\mapsto I_{d+1}+u.e_{ij}\label{alphaij}\end{gather}where $I_{d+1}+u.e_{ij}$ is the matrix $(u_{rs})$ with $u_{rr}=1$ (all $r$), with $u_{ij}=u$ and with $u_{rs}=0$ for all other pairs $(r,s)$. For the root $\alpha=\epsilon_i-\epsilon_j\in \Phi$ and $r\in\mathbb{R}$ let $$U_{\alpha,r}=\widetilde{\alpha}_{ij}(\{u\in K;\,\omega(u)\ge r\})\subset G.$$ For $x\in X_*(T)\otimes{\mathbb{R}}$ let$$U_x=\mbox{the subgroup of}\,\,G\,\,\mbox{generated by all}\,\,U_{\alpha,-\alpha(x)}\,\,\mbox{for}\,\,\alpha\in\Phi.$$Let $W$ be the subgroup of permutation matrices in $G$ and let $N=T\rtimes W$ be the normalizer of $T$ in $G$.

Let now $g\in G$ and $Z\in F^0_{A}$ such that also $g.Z\in F^0_{A}$. We claim that $g:\underline{M}_{\dot{K}}|_{\mathfrak{U}_{Z}}\cong \underline{M}_{\dot{K}}|_{\mathfrak{U}_{g.Z}}$ restricts to an isomorphism$$g:{\mathcal M}^0_{{\mathcal O}_{\dot{K}}}|_{\mathfrak{U}_{Z}}\cong{\mathcal M}^0_{{\mathcal O}_{\dot{K}}}|_{\mathfrak{U}_{g.Z}}.$$We have a canonical projection from $X_*(T)\otimes{\mathbb{R}}$ to the standard apartment in the Bruhat Tits building of ${\rm PGL}\sb {d+1}/K$ (see \cite{schtei}). Suppose that $x\in X_*(T)\otimes{\mathbb{R}}$ projects to the vertex corresponding to $Z\in F^0_{A}$. (In the above mentioned correspondence between $F^0_{A}$ and vertices in the standard apartment.) By the Bruhat decomposition, there exist $h_x\in U_x$, $h_{gx}\in U_{gx}$ and $n\in N$ such that $g=h_{gx}nh_{x}$. Therefore we may split up our task into the following cases (i)-(iii):\\
(i) $g\in T$,\\(ii) $g\in W$,\\ (iii) $x=gx$ and $g\in U_x$ for some $x\in X_*(T)\otimes{\mathbb{R}}$.

(i) Suppose $g\in T$. We claim that in this case $g$ even respects weight spaces: we prove that $g$ induces for any $\mu\in X^*(T)$ with $M_{\mu}\ne0$ an isomorphism$$g:\dot{\pi}^{(d+1)\overline{\mu}(Z)}M_{\mu}^0\otimes_{{\mathcal O}_K}{\mathcal O}_{\dot{K}}\cong\dot{\pi}^{(d+1)\overline{\mu}(g.Z)}M_{\mu}^0\otimes_{{\mathcal O}_K}{\mathcal O}_{\dot{K}}.$$Indeed, $g$ induces an isomorphism $$g:M_{\mu}^0\cong\pi^{\omega(\mu(g))}M_{\mu}^0.$$Thus, according to our conventions regarding the action of $G$ on $M\otimes_K\dot{K}$, it induces an isomorphism$$g:M_{\mu}^0\otimes_{{\mathcal O}_K}{\mathcal O}_{\dot{K}}\cong\dot{\pi}^{(d+1)\omega(\mu(g))-|M|\omega({\rm{det}}(g))}M_{\mu}^0\otimes_{{\mathcal O}_K}{\mathcal O}_{\dot{K}}.$$But $$(d+1)\omega(\mu(g))-|M|\omega({\rm{det}}(g))=(d+1)(-\omega(\overline{\mu}(g)))=(d+1)(\overline{\mu}(g.Z)-\overline{\mu}(Z))$$and the claim follows. 

(ii) Now $g\in W$. The isomorphisms $g:M_{\mu}\cong M_{g.\mu}$ restrict to isomorphisms $g:M_{\mu}^0\cong M_{g\mu}^{0}$ since $M^0\subset M$ is stable under ${\rm GL}\sb {d+1}({\mathcal O}_K)$. On the other hand $\mu(Z)=(g.\mu)(g.Z)$ and hence $\overline{\mu}(Z)=\overline{(g.\mu)}(g.Z)$ for $\mu\in X^*(T)$. It follows that $g$ induces isomorphisms$$\dot{\pi}^{(d+1)\overline{\mu}(Z)}M_{\mu}^{0}\otimes_{{\mathcal O}_K}{\mathcal O}_{\dot{K}}\cong\dot{\pi}^{(d+1)\overline{(g.\mu)}(g.Z)}M_{g.\mu}^{0}\otimes_{{\mathcal O}_K}{\mathcal O}_{\dot{K}}$$for any $\mu\in X^*(T)$ and we are done for such $g$.

(iii) Now consider the case $x=g.x$ and $g\in U_x$ for some $x\in X_*(T)\otimes{\mathbb{R}}$. Then also $Z=g.Z$ and $\overline{\mu}(x)=\overline{\mu}(Z)$. We may assume $g\in U_{\alpha,-\alpha(x)}$ for some $\alpha=\epsilon_i-\epsilon_j\in\Phi$. Thus $g=\widetilde{\alpha}_{ij}(u)$ for some $u\in K$ with $\omega(u)\ge-\alpha(x)$. It suffices to show that the automorphism $g$ of $M$ induces an automorphism$$g:\bigoplus_{\mu\in X^*(T)}\dot{\pi}^{(d+1)\overline{\mu}(x)}M_{\mu}^0\cong\bigoplus_{\mu\in X^*(T)}\dot{\pi}^{(d+1)\overline{\mu}(x)}M_{\mu}^0.$$Now $\omega(u)\ge-\alpha(x)$ implies $\overline{\mu+n\alpha}(x)\le\overline{\mu}(x)+n\omega(u)$ for all $\mu\in X^*(T)$, all $n\in\mathbb{N}_0$. Therefore it is enough to prove
\begin{gather}\widetilde{\alpha}_{ij}(u)(m)\subset\sum_{n\ge0}u^n.M_{\mu+n(\epsilon_i-\epsilon_j)}^0\label{liearg}\quad\quad (m\in M_{\mu}^0).\end{gather}To see this define $X_{\alpha}=(d\widetilde{\alpha}_{ij})(1)\in\mbox{Lie}({\rm GL}\sb {d+1}/\mathbb{Z})$ for $\alpha=\epsilon_i-\epsilon_j$ and then $$X_{\alpha,n}=\frac{X_{\alpha}^n}{n!}\in\mbox{Dist}({\rm GL}\sb {d+1}/\mathbb{Z})\quad\mbox{for}\quad n\ge0$$(compare \cite{jan} II.1.11 and 1.12). By \cite{jan} II.1.19 we have$$X_{\alpha,n}M_{\mu}\subset M_{\mu+n\alpha}\quad\mbox{and}\quad\widetilde{\alpha}_{ij}(u)(m)=\sum_{n\ge0}u^nX_{\alpha,n}(m).$$Since $X_{\alpha,n}$ is defined over $\mathbb{Z}$ we in turn have $X_{\alpha,n}M_{\mu}^0\subset M_{\mu+n\alpha}^0$ and formula (\ref{liearg}) follows.

The above claim is established. It follows that on the dense open formal subscheme ${\bigcup_{Z\in F^0}\mathfrak{U}_{Z}}$ of $\mathfrak{X}$ (the union is disjoint) there is a unique $G$-stable subsheaf $${\mathcal M}^0_{{\mathcal O}_{\dot{K}}}|_{\bigcup_{Z\in F^0}\mathfrak{U}_{Z}}\subset\underline{M}_{\dot{K}}|_{\bigcup_{Z\in F^0}\mathfrak{U}_Z}$$whose restriction to ${\mathfrak{U}_{Z}}$ for $Z\in F^0_{A}$ is ${\mathcal M}^0_{{\mathcal O}_{\dot{K}}}|_{\mathfrak{U}_{Z}}$ as defined by (\ref{numdef}). We define ${\mathcal M}^0_{{\mathcal O}_{\dot{K}}}$ on all of $\mathfrak{X}$ as the maximal ${\mathcal O}_{\dot{K}}$-module subsheaf of $\underline{M}_{\dot{K}}$ restricting to ${\mathcal M}^0_{{\mathcal O}_{\dot{K}}}|_{\cup_{Z\in F^0}\mathfrak{U}_{Z}}$.\hfill$\Box$\\

We now wish to compute the reduction modulo $\dot{\pi}$ of ${\mathcal M}^0_{{\mathcal O}_{\dot{K}}}\otimes_{{\mathcal O}_K}{\mathcal O}_{{\mathfrak{X}}}$ in terms of our sheaves ${\widetilde{\mathcal M}}[{\mathcal O}_Y]$ living on the central irreducible component $Y$.

For open formal subschemes ${\mathfrak{U}}$ of ${\mathfrak{X}}$ we write ${\mathcal O}_{\dot{{\mathfrak{U}}}}={\mathcal O}_{\mathfrak{U}}\otimes_{{\mathcal O}_K}{\mathcal O}_{\dot{K}}$. Recall that in section \ref{rera} we defined the open subscheme $Y'$ of $Y$ with $U(k)Y'=Y$, defined the $T(k)$-equivariant sheaf $\widetilde{\mathcal M}[{\mathcal O}_{Y'}]$ on $Y'$ with a $T(k)$-equivariant $X^*(T)$-indexed grading and extended it to a ${\rm GL}_{d+1}(k)$-equivariant sheaf $\widetilde{\mathcal M}[{\mathcal O}_{Y}]$ with a $U(k)$-equivariant $X^*(T)$-indexed filtration on $Y$. We now perform a similar construction in our present global setting. Here the role of $Y'$ is played by ${\mathfrak Y}$: by definition, ${\mathfrak Y}$ is the open formal subscheme of ${\mathfrak X}$ such that for the open subscheme ${\mathfrak Y}\otimes k$ of ${\mathfrak X}\otimes k$ we have$${\mathfrak X}\otimes k-{\mathfrak Y}\otimes k=\bigcup_{Z\in F^0-F^0_{A}}Z.$$We have $U(K).{\mathfrak Y}={\mathfrak X}$. Moreover observe $Y'={\mathfrak Y}\cap Y$.

For $Z\in F^0_{A}$ let ${\mathcal J}_{Z}\subset {\mathcal O}_{\mathfrak Y}$ be the ideal defining the closed subscheme $Z\cap {\mathfrak Y}$ inside ${\mathfrak Y}$. Note that ${\mathcal J}_{Z}$ is invertible inside ${\mathcal O}_{\mathfrak Y}\otimes_{{\mathcal O}_K}K$: indeed, small open formal subschemes of ${\mathfrak Y}$ admit open embeddings into the $\pi$-adic completion of $\spec({\mathcal O}_K[X_0,\ldots,X_d]/(X_0\ldots X_d-\pi))$, and for an appropriate numbering of $X_0,\ldots,X_d$ the element $X_0$ is a generator of ${\mathcal J}_{Z}$; in $K[X_0,\ldots,X_d]/(X_0\ldots X_d-\pi)$ its inverse is $\pi^{-1}X_1\ldots X_d$. Thus we may speak of negative integral powers of ${\mathcal J}_{Z}$ as ${\mathcal O}_{\mathfrak Y}$-submodules of ${\mathcal O}_{\mathfrak Y}\otimes_{{\mathcal O}_K}K$. Also note that on small open formal subschemes of ${\mathfrak Y}$ we have ${\mathcal J}_{Z}={\mathcal O}_{\mathfrak Y}$ for almost all $Z$, therefore the following infinite products of ${\mathcal O}_{{\mathfrak Y}}$-submodules in ${\mathcal O}_{\mathfrak Y}\otimes_{{\mathcal O}_K}{\dot{K}}$ make sense. For any $\mu\in X^*(T)$ we define the sheaf\begin{gather}({\mathcal O}_{\dot{\mathfrak Y}})^{\overline{\mu}}=\sum_{s=0}^d\dot{\pi}^s\prod_{Z\in F^0_{A}}{\mathcal J}_{Z}^{\lceil\overline{\mu}(Z)-\frac{s}{d+1}\rceil},\label{iii}\end{gather}on ${\mathfrak Y}$, the ${\mathcal O}_{\dot{\mathfrak Y}}$-submodule of ${\mathcal O}_{{\mathfrak{Y}}}\otimes_{{\mathcal O}_K}{\dot{K}}$ generated by the submodules $\dot{\pi}^s\prod_{Z\in F^0_{A}}{\mathcal J}_{Z}^{\lceil\overline{\mu}(Z)-\frac{s}{d+1}\rceil}$ for $s=0,\ldots,d$. Let $({\mathcal O}_{\dot{\mathfrak X}})^{\overline{\mu}}$ be the unique ${U}(K)$-equivariant ${\mathcal O}_{\dot{\mathfrak X}}$-module subsheaf of ${\mathcal O}_{{\mathfrak{X}}}\otimes_{{\mathcal O}_K}{\dot{K}}$ (with its ${U}(K)$-action induced by that of $G$ on ${\mathfrak{X}}$) whose restriction to ${\mathfrak Y}$ is $({\mathcal O}_{\dot{\mathfrak Y}})^{\overline{\mu}}$. To describe the reduction of $({\mathcal O}_{\dot{\mathfrak X}})^{\overline{\mu}}$ we need to parametrize the set ${\mathcal V}$ in terms of the action of $U(k)$ on it. For $\sigma\subsetneq\{0,\ldots,d\}$ let$${U}_{\sigma}=\{(a_{ij})_{0\le i,j\le d}\in{U}\quad|\quad a_{ij}=0\,\mbox{if}\,i\ne j\,\mbox{and}\,[j\notin\sigma\,\mbox{or}\,\{i,j\}\subset \sigma]\}.$$Let
$${\mathcal N}=\{(\sigma,u)\quad|\quad\emptyset\ne\sigma\subsetneq\{0,\ldots,d\}, u\in {U}_{\sigma}(k)\}.$$We have a bijection (see \cite{holdis})$${\mathcal N}\cong{\mathcal V},\quad\quad (\sigma,u)\mapsto u.V_{{\sigma}}$$and the set of orbits of ${U}(k)$ acting on the set ${\mathcal V}$ is in bijection with the set of all $\sigma$ with $\emptyset\ne\sigma\subsetneq\{0,\ldots,d\}$.

We will need the sheaves $({\mathcal O}_{\dot{\mathfrak X}})^{\overline{\mu}}$ only for those $\mu$ with $M_{\mu}\ne0$. For such $\mu$ consider the partition of $F^0_A$, indexed by the $t\in\{0,1,\ldots,d\}$, into the subsets$$F_A^0(t)=\{Z\in F_A^0 \quad | \quad \overline{\mu}(Z)-\frac{t}{d+1}\in\mathbb{Z}\}.$$It provides the partition of $F^0$ into the subsets$$F^0(t)=U(K).F_A^0(t).$$Since $M$ is irreducible, all $\mu$ with $M_{\mu}\ne0$ differ by linear combinations of elements of $\Phi$ (see \cite{jan} II.2.2). For each such $\mu$, if we write $\overline{\mu}=\sum_{j=0}^d\overline{a}_j(\mu)\epsilon_j$ (cf. formula (\ref{aquer})), we have \begin{gather}\overline{a}_j(\mu)-\frac{|M|}{d+1}\in\mathbb{Z}\label{part}\end{gather}for all $0\le j\le d$. It follows that $F_A^0(t)$ and hence $F^0(t)$ does not depend on $\mu$ and moreover that $F^0(t)$ is non-empty if ond only if $t\equiv n|M|$ modulo $(d+1)$ for some $n\in\mathbb{Z}$, or in other words: we have defined a partition of $F^0$ indexed by the multiples of (the class of) $|M|$ in $\mathbb{Z}/(d+1)$. This partition is stable for the action of ${\rm SL}_{d+1}(K)$ on $F^0$ (this again follows from equation (\ref{part})) while the action of the full group $G$ on $F^0$ can be used to cycle through the parts of this partition. Endow $$\overline{\mathfrak{X}(t)}=\bigcup_{Z\in F^0(t)}Z$$with its structure of reduced closed subscheme of ${\mathfrak{X}}\otimes k$. 

\begin{lem}\label{lbdlext} We have natural isomorphisms\begin{gather}({\mathcal O}_{\dot{\mathfrak X}})^{\overline{\mu}}\otimes_{{\mathcal O}_{\dot{K}}}k\cong\bigoplus_{t\in\{0,\ldots,d\}}   \frac{({\mathcal O}_{\dot{\mathfrak X}})^{\overline{\mu}}\otimes_{{\mathcal O}_{\dot{\mathfrak{X}}}}{\mathcal O}_{\overline{\mathfrak{X}(t)}}}{{\mathcal O}_{\overline{\mathfrak{X}(t)}}\mbox{-torsion}},\label{lbdlzerldeninu}\end{gather}\begin{gather}{\mathcal L}_{Y}(\sum_{\emptyset\ne\sigma\subsetneq\{0,\ldots,d\}}-\lceil\sum_{j\in\sigma}\overline{a}_j(\mu)\rceil\sum_{u\in{U}_{\sigma}(k)}u.V_{\sigma})\cong \frac{({\mathcal O}_{\dot{\mathfrak X}})^{\overline{\mu}}\otimes_{{\mathcal O}_{\dot{\mathfrak{X}}}}{\mathcal O}_{Y}}{{\mathcal O}_{Y}\mbox{-torsion}}.\label{lbdldeninu}\end{gather}
\end{lem}

{\sc Proof:} Let ${\mathfrak Y}(t)$ denote the maximal open formal subscheme of ${\mathfrak Y}$ such that the open subscheme ${\mathfrak Y}(t)\otimes k$ of ${\mathfrak Y}\otimes k$ is contained in $\cup_{Z\in F_A^0(t)}(Z\cap {\mathfrak Y})$. Let $$\overline{{\mathfrak Y}(t)}=\bigcup_{Z\in F^0(t)}(Z\cap{\mathfrak Y})$$with its reduced structure, or equivalently: $\overline{{\mathfrak Y}(t)}$ is the schematic closure of ${\mathfrak Y}(t)\otimes k$ in ${\mathfrak Y}\otimes k$. By formula (\ref{iii}) the restriction of $({\mathcal O}_{\dot{\mathfrak Y}})^{\overline{\mu}}$ to ${\mathfrak Y}(t)$ is the line bundle$$\dot{\pi}^{t}\prod_{Z\in F^0_{A}}{\mathcal J}_{Z}^{\lceil\overline{\mu}(Z)-\frac{t}{d+1}\rceil}|_{{\mathfrak Y}(t)}=\dot{\pi}^{t}\prod_{Z\in F^0_{A}(t)}{\mathcal J}_{Z}^{\overline{\mu}(Z)-\frac{t}{d+1}}|_{{\mathfrak Y}(t)}$$on ${\mathfrak Y}(t)$ (all ${\mathcal J}_{Z}|_{{\mathfrak Y}(t)}$ for $Z\in F^0_A-F^0_A(t)$ are trivial). We obtain: the reduction $({\mathcal O}_{\dot{\mathfrak Y}})^{\overline{\mu}}\otimes_{{\mathcal O}_{\dot{K}}}k$ of $({\mathcal O}_{\dot{\mathfrak Y}})^{\overline{\mu}}$ decomposes into a direct sum, indexed by the set $\{0,\ldots,d\}$, whose direct summand for $t\in\{0,\ldots,d\}$ is the image of the map$$\dot{\pi}^{t}\prod_{Z\in F^0_{A}}{\mathcal J}_{Z}^{\lceil\overline{\mu}(Z)-\frac{t}{d+1}\rceil}\longrightarrow({\mathcal O}_{\dot{\mathfrak Y}})^{\overline{\mu}}\longrightarrow{({\mathcal O}_{\dot{\mathfrak Y}})^{\overline{\mu}}}\otimes_{{\mathcal O}_{\dot{\mathfrak{Y}}}}{\mathcal{O}}_{\overline{{\mathfrak Y}(t)}}.$$This is a line bundle on $\overline{{\mathfrak Y}(t)}$ and maps isomorphically to the quotient of ${({\mathcal O}_{\dot{\mathfrak Y}})^{\overline{\mu}}}\otimes_{{\mathcal O}_{\dot{\mathfrak{Y}}}}{\mathcal{O}}_{\overline{{\mathfrak Y}(t)}}$ divided by its ${\mathcal O}_{\overline{{\mathfrak Y}(t)}}$-torsion. Thus\begin{gather}({\mathcal O}_{\dot{\mathfrak Y}})^{\overline{\mu}}\otimes_{{\mathcal O}_{\dot{K}}}k\cong\bigoplus_{t\in\{0,\ldots,d\}}\frac{({\mathcal O}_{\dot{\mathfrak Y}})^{\overline{\mu}}\otimes_{{\mathcal O}_{\dot{\mathfrak{Y}}}}{\mathcal{O}}_{\overline{{\mathfrak Y}(t)}}}{{\mathcal{O}}_{\overline{{\mathfrak Y}(t)}}\mbox{-torsion}}\label{deltaninu}\end{gather}and the direct summand for $t\in\{0,\ldots,d\}$ is an invertible ${\mathcal{O}}_{\overline{{\mathfrak Y}(t)}}$-module. Hence formula (\ref{lbdlzerldeninu}) by  ${U}(K)$-equivariance. We also see $$\frac{({\mathcal O}_{\dot{\mathfrak X}})^{\overline{\mu}}\otimes_{{\mathcal O}_{\dot{\mathfrak{X}}}}{\mathcal O}_{\overline{\mathfrak{X}}(0)}}{{\mathcal O}_{\overline{\mathfrak{X}}(0)}\mbox{-torsion}}\otimes_{{\mathcal O}_{\overline{\mathfrak{X}}(0)}}{\mathcal O}_{Y}=\frac{({\mathcal O}_{\dot{\mathfrak X}})^{\overline{\mu}}\otimes_{{\mathcal O}_{\dot{\mathfrak{X}}}}{\mathcal O}_{Y}}{{\mathcal O}_{Y}\mbox{-torsion}}$$and that this is the unique ${U}(k)$-equivariant subsheaf of the constant sheaf $k(Y)$ on $Y$ whose restriction to $Y'=Y\cap\mathfrak{Y}$ is $$\frac{({\mathcal O}_{\dot{\mathfrak Y}})^{\overline{\mu}}\otimes_{{\mathcal O}_{\dot{\mathfrak{Y}}}}{\mathcal{O}}_{Y'}}{{\mathcal{O}}_{Y'}\mbox{-torsion}}$$(for the uniqueness note that ${U}(k).Y'=Y$). If we define the divisor $D$ on $Y'$ by requiring$$\frac{({\mathcal O}_{\dot{\mathfrak Y}})^{\overline{\mu}}\otimes_{{\mathcal O}_{\dot{\mathfrak{Y}}}}{\mathcal{O}}_{Y'}}{{\mathcal{O}}_{Y'}\mbox{-torsion}}={\mathcal L}_{Y'}(D)$$(as subsheaves of the constant sheaf $k(Y)$ on $Y'$), then by ${U}(k)$-equivariance of its both sides and ${U}(k).Y'=Y$, to prove formula (\ref{lbdldeninu}) we only need to prove the identity of divisors$$D=\sum_{\emptyset\ne\sigma\subsetneq\{0,\ldots,d\}}-\lceil\sum_{j\in\sigma}\overline{a}_j(\mu)\rceil V_{\sigma}|_{Y'}$$on $Y'$. To see this note that for $\emptyset\ne \sigma\subsetneq\{0,\ldots,d\}$ we have $V_{\sigma}=Z_{{\sigma}}\cap Y$ on $Y\in F^0_A$; here we write $Z_{\sigma}=t_{\sigma}Y\in F^0_A$ with $t_{\sigma}\in T\subset G$ defined as $t_{\sigma}=\diag(t_{\sigma,0},\ldots,t_{\sigma,d})$ with $t_{\sigma,j}=1$ if $j\notin\sigma$ and $t_{\sigma,j}=\pi$ if $j\in\sigma$. Now we only need to see that for $\emptyset\ne\sigma\subsetneq\{0,\ldots,d\}$ the prime divisor $V_{\sigma}=Z_{{\sigma}}\cap Y$ occurs in $D$ with multiplicity $$-\lceil\overline{\mu}(Z_{\sigma})\rceil=-\lceil\sum_{j\in\sigma}\overline{a}_j(\mu)\rceil.$$But our discussion shows that $\frac{({\mathcal O}_{\dot{\mathfrak Y}})^{\overline{\mu}}\otimes_{{\mathcal O}_{\dot{\mathfrak{Y}}}}{\mathcal{O}}_{Y'}}{{\mathcal{O}}_{Y'}\mbox{-torsion}}$ can be identified with the image of the map$$\prod_{\emptyset\ne\sigma\subsetneq\{0,\ldots,d\}}{\mathcal J}_{Z_{\sigma}}^{\lceil\overline{\mu}(Z_{\sigma})\rceil}\longrightarrow({\mathcal O}_{\dot{\mathfrak Y}})^{\overline{\mu}}\longrightarrow{({\mathcal O}_{\dot{\mathfrak Y}})^{\overline{\mu}}}\otimes_{{\mathcal O}_{\dot{\mathfrak{Y}}}}{\mathcal{O}}_{Y'}$$and we can read off the correct multiplicity.\hfill$\Box$

\begin{satz}\label{strumn} Let $\iota_Y:Y\to{\mathfrak X}$ denote the closed embedding. We have natural isomorphisms$$({\mathcal M}^0_{{\mathcal O}_{\dot{K}}}\otimes_{{\mathcal O}_K}{\mathcal O}_{{\mathfrak{X}}})\otimes_{{\mathcal O}_{\dot{K}}}k\cong\bigoplus_{t\in\{0,\ldots,d\}}\frac{({\mathcal M}^0_{{\mathcal O}_{\dot{K}}}\otimes_{{\mathcal O}_K}{\mathcal O}_{{\mathfrak{X}}})\otimes_{{\mathcal O}_{\dot{\mathfrak{X}}}}{\mathcal O}_{\overline{\mathfrak{X}(t)}}}{{\mathcal O}_{\overline{\mathfrak{X}(t)}}\mbox{-torsion}}$$\begin{gather}\frac{({\mathcal M}^0_{{\mathcal O}_{\dot{K}}}\otimes_{{\mathcal O}_K}{\mathcal O}_{{\mathfrak{X}}})\otimes_{{\mathcal O}_{\dot{\mathfrak{X}}}}{\mathcal O}_Y}{{\mathcal O}_Y\mbox{-torsion}}\cong \iota_{Y,*}({\widetilde{\mathcal M}}[{\mathcal O}_Y]).\label{kofo}\end{gather}
\end{satz}

{\sc Proof:} To prove Theorem \ref{strumn} it suffices by $U(K)$-equivariance (resp. by $U({\mathcal O}_K)$-equivariance for the isomorphism (\ref{kofo})) to prove the statements on the sheaves restricted to ${\mathfrak Y}$ (resp. to $Y'$ for the isomorphism (\ref{kofo})). There, by construction, ${\mathcal M}^0_{{\mathcal O}_{\dot{K}}}\otimes_{{\mathcal O}_K}{\mathcal O}_{{\mathfrak{X}}}$ decomposes into a direct sum indexed by the weights $\mu$ of $M$. A small computation in local coordinates shows that formula (\ref{numdef}) implies that the summand for $\mu$ is of the form $M_{\mu}^0\otimes_{{\mathcal O}_K}({\mathcal O}_{\dot{\mathfrak X}})^{\overline{\mu}}$ so that Lemma \ref{lbdlext} proves the first isomorphism. The isomorphism (\ref{kofo}) now follows from formula (\ref{ystrdf}) (in view of Lemma \ref{lbdlext}).\hfill$\Box$\\

{\it Remark:} If $|M|\in(d+1).\mathbb{Z}$, or equivalently if $\overline{\mu}\in X^*(\overline{T})$ for all $\mu$ with $M_{\mu}\ne0$, then $F^0=F^0(0)$ and $({\mathcal O}_{\dot{\mathfrak X}})^{\overline{\mu}}$ is a line bundle on ${\mathfrak X}$ for each such $\mu$, and ${\mathcal M}^0_{{\mathcal O}_{\dot{K}}}\otimes_{{\mathcal O}_K}{\mathcal O}_{{\mathfrak{X}}}$ is a locally free ${\mathcal O}_{{\mathfrak{X}}}$-module sheaf on ${\mathfrak{X}}$.

\section{Coherent cohomology via logarithmic differential forms}
\label{hodese}

Let ${\mathfrak S}$ be a strictly semistable formal ${\mathcal O}_K$-scheme. Endow ${\mathfrak S}$ and $\spf({\mathcal O}_K)$ with the log structure defined by the respective special fibre and let $\Omega^{\bullet}_{{\mathfrak S}}$ denote the logarithmic de Rham complex for the log smooth morphism of formal log schemes ${\mathfrak S}\to\spf({\mathcal O}_K)$. Let $\Omega^{\bullet}_S$ denote the push forward to ${\mathfrak S}$ of the de Rham complex on $S={\mathfrak S}\otimes_{{\mathcal O}_K}K$ (a $K$-rigid space). Let $T$ be an irreducible component of the special fibre ${\mathfrak S}\otimes k$ of ${\mathfrak S}$, let $T^0$ denote the maximal open subscheme of ${\mathfrak S}\otimes k$ which is contained in $T$. Then $T-T^0$ is a normal crossings divisor on the smooth $k$-scheme $T$. Let $\Omega_T^{\bullet}$ denote the de Rham complex on $T$ with logarithmic poles along $T-T^0$.
 
\begin{lem}\label{logderacom} There are canonical isomorphism of sheaf complexes$$\Omega^{\bullet}_S\cong\Omega^{\bullet}_{{\mathfrak S}}\otimes_{{\mathcal O}_K}K,\quad\quad\quad\Omega_T^{\bullet}\cong\Omega^{\bullet}_{{\mathfrak S}}\otimes_{{\mathcal O}_{{\mathfrak S}}}{\mathcal O}_T.$$
\end{lem}

{\sc Proof:} The first isomorphism is clear. To prove the second one let $T'$ (resp. $\spec(k)'$) denote the log scheme whose underlying scheme is $T$ (resp. $\spec(k)$) and whose log structure is the pull back of that of ${\mathfrak S}$ (resp. that of $\spf({\mathcal O}_K)$). In other words, $T'\to {\mathfrak S}$ and $\spec(k)'\to \spf({\mathcal O}_K)$ are {\it exact} closed immersions of log schemes. Then $T'$ is a log scheme over the base $\spec(k)'$ (in general not log smooth). Let $\Omega_{T'}^{1}$ be the logarithmic differential module of $T'\to\spec(k)'$. We have a morphism of log schemes $T'\to T$. By functoriality we get natural morphisms of sheaves$$\Omega^{1}_{T}\longrightarrow\Omega^{1}_{T'}\longleftarrow\Omega^{1}_{{\mathfrak S}}\otimes_{{\mathcal O}_{{\mathfrak S}}}{\mathcal O}_T.$$We claim that both are isomorphisms. To see this we may assume that ${\mathfrak S}$ is the formal $\pi$-adic completion of $\spec({\mathcal O}_K[X_0,\ldots,X_d]/(X_0\cdots X_s-\pi))$ for some $0\le s\le d$ and that the kernel of ${\mathcal O}_{{\mathfrak S}}\to{\mathcal O}_T$ is generated by $X_0$. Then these sheaves are canonically identified with the free ${\mathcal O}_T$-module with basis $\{\dlog(X_1),\ldots,\dlog(X_s),{\rm d}(X_{s+1}),\ldots,{\rm d}(X_d)\}$. The lemma follows.\hfill$\Box$\\

\begin{lem}\label{comalg} Let $A$ be a discrete valuation ring with uniformizer $\lambda\in A$, residue field $k$ and fraction field $F$. Let $M$ be a $\lambda$-torsion free $A$-module.\\(a) The $A$-module $$M'=\lim_{\leftarrow\atop n}M/\lambda^nM$$ is $\lambda$-torsion free. For each $r\ge1$ the map $M/\lambda^rM\to M'/\lambda^rM'$ induced by the natural map $M\to M'$ is bijective; in particular we have $$M'=\lim_{\leftarrow\atop n}M'/\lambda^nM'.$$(b) Let $\widetilde{N}$ be a sub vector space of $M\otimes_AF$ and let $N=M\cap\widetilde{N}$ (intersection inside $M\otimes_AF$). The map $N\otimes_Ak\to M\otimes_Ak$ induced by the natural map $N\to M$ is injective. 
\end{lem}

{\sc Proof:} (a) Suppose we are given $(m_n)_n\in M'$ and $s\ge1$ such that $\lambda^s(m_n)_n=0$ in $M'$. Let $n\ge1$ and choose $x\in M$ such that $\overline{x}=m_{n+s}\in M/\lambda^{n+s}M$ (where $\overline{x}$ denotes the image of $x$ in $M/\lambda^{n+s}M$). Then $\lambda^sm_{n+s}=0$ in $M/\lambda^{n+s}M$ implies $\lambda^sx\in\lambda^{n+s}M$, hence $x\in \lambda^{n}M$ since $M$ is $\lambda$-torsion free, hence $m_n=0$ in $M/\lambda^{n}M$ and we see that $M'$ is $\lambda$-torsion free. Next let $r\ge1$ and suppose we are given $(m_n)_n\in M'$. Let $m\in M$ be an arbitrary lift of $m_r\in M/\lambda^r M$. We find an element $(b_n)_n\in M'$ such that $\lambda^r b_n=\overline{m}-m_n\in M/\lambda^nM$ for all $n$ (here $\overline{m}$ denotes the class of $m$). Indeed, we know $\overline{m}-m_{n+r}\in \lambda^r M/\lambda^{n+r}M$. Choose $b'_n\in M$ with $\lambda^r b'_n=\overline{m}-m_{n+r}\in M/\lambda^{n+r}M$ and let $b_n$ be the image of $b'_n$ in $M/\lambda^{n}M$. Then $(b_n)_n\in M'$ because $(m_n)_n\in M'$ implies $\lambda^r(b'_{n+r}-b'_n)\in\lambda^{n+r}M$, hence $b'_{n+r}-b'_n\in\lambda^{n}M$ since $M$ is $\lambda$-torsion free. Now $\lambda^r((b_n)_n)=(\overline{m}-m_n)_n$ in $M'$, thus $m\in M$ and $(m_n)_n\in M'$ map to the same element in $M'/\lambda^r M'$. We have shown that $M/\lambda^rM\to M'/\lambda^rM'$ is surjective; the injectivity is clear.\\(b) This is very easy.\hfill$\Box$\\

In the sequel, for sheaves ${\mathcal G}$ on $X$ we write ${\mathcal G}$ also for the push forward sheaf on $\mathfrak{X}$ under the specialization map $sp:X\to \mathfrak{X}$; we use tacitly and repeatedly Kiehl's result \cite{kiaub} that if ${\mathcal G}$ is coherent we have ${\mathbb R}^tsp_*{\mathcal G}=0$ for all $t>0$. 

Denote by $\Omega^{\bullet}_{\mathfrak{X}}$ the logarithmic de Rham complex of the log smooth morphism of formal log schemes ${\mathfrak{X}}\to\spf({\mathcal O}_K)$, where we give the source and the target the respective log structures defined by the special fibres. Note that by Lemma \ref{logderacom} we have canonical identifications $$\Omega^{\bullet}_{\mathfrak{X}}\otimes_{{\mathcal O}_K}K=\Omega^{\bullet}_X,\quad\quad\quad\quad\Omega^{\bullet}_{\mathfrak{X}}\otimes_{{\mathcal O}_{\mathfrak{X}}}{{\mathcal O}_Y}=\Omega^{\bullet}_Y.$$

Recall that we view $X$ as a subspace of ${\mathbb P}_K^d$. For $0\le s\le d$ let ${\mathcal Log}^s$ be the $K$-vector subspace of $\Omega_X^s(X)$ generated by logarithmic differential $s$-forms on ${\mathbb P}_K^d$ with logarithmic poles along $K$-rational hyperplanes. In other words, ${\mathcal Log}^s$ is generated by $s$-forms $\eta$ of the type\begin{gather}\eta=\dlog(y_1)\wedge\ldots\wedge\dlog(y_s)\label{frofo}\end{gather}for which there exists a suitable (adapted to $\eta$) choice of projective coordinate system $\theta_0,\ldots,\theta_d$ on ${\mathbb P}_K^d$ (i.e. a suitable (adapted to $\eta$) isomorphism of $K$-varieties ${\mathbb P}_K^d\cong\proj(K[\theta_0,\ldots,\theta_d])$) such that $y_j=\theta_j/\theta_0\in {\mathcal O}_X^{\times}(X)$ for all $1\le j\le s$. Clearly ${\mathcal Log}^s$ is a $G$-stable subspace of $\Omega_X^s(X)$. 

Let ${\mathcal M}^0_{{\mathcal O}_{\dot{K}}}$ be the $G$-equivariant integral structure in the constant sheaf $\underline{M}_{\dot{K}}$ defined in section \ref{secints}. For an open quasi-compact subscheme $U$ of ${\mathfrak X}\otimes k$ we have $M\otimes_K\dot{K}\otimes_K\Omega_X^s(U)=({\mathcal M}^0_{{\mathcal O}_{\dot{K}}}\otimes_{{\mathcal O}_K}\Omega^{s}_{{\mathfrak{X}}})(U)\otimes K$, hence for such $U$ we may view $M\otimes_K\dot{K}\otimes_K\Omega_X^s(X)$ and consequently also $M\otimes_K\dot{K}\otimes_K{\mathcal Log}^s$ as being contained in $({\mathcal M}^0_{{\mathcal O}_{\dot{K}}}\otimes_{{\mathcal O}_K}\Omega^{s}_{{\mathfrak{X}}})(U)\otimes K$. Therefore we may define$${\mathcal Log}_{alg}^s({\mathcal M}^0_{{\mathcal O}_{\dot{K}}})(U)=M\otimes_K\dot{K}\otimes_K{\mathcal Log}^s\quad\bigcap\quad({\mathcal M}^0_{{\mathcal O}_{\dot{K}}}\otimes_{{\mathcal O}_K}\Omega^{s}_{{\mathfrak{X}}})(U),$$the intersection taking place inside $({\mathcal M}^0_{{\mathcal O}_{\dot{K}}}\otimes_{{\mathcal O}_K}\Omega^{s}_{{\mathfrak{X}}})(U)\otimes K$. Since the restriction maps of the sheaf ${\mathcal M}^0_{{\mathcal O}_{\dot{K}}}\otimes_{{\mathcal O}_K}\Omega^{s}_{{\mathfrak{X}}}$ are injective we have thus defined a $G$-stable subsheaf ${\mathcal Log}_{alg}^s({\mathcal M}^0_{{\mathcal O}_{\dot{K}}})$ of ${\mathcal M}^0_{{\mathcal O}_{\dot{K}}}\otimes_{{\mathcal O}_K}\Omega^{s}_{{\mathfrak{X}}}$. For open $U\subset{\mathfrak X}\otimes k$ we further let $${\mathcal Log}^s({\mathcal M}^0_{{\mathcal O}_{\dot{K}}})(U)=\lim_{\leftarrow\atop n}\frac{{\mathcal Log}_{alg}^s({\mathcal M}^0_{{\mathcal O}_{\dot{K}}})(U)}{\dot{\pi}^n{\mathcal Log}_{alg}^s({\mathcal M}^0_{{\mathcal O}_{\dot{K}}})(U)}.$$This defines a sheaf ${\mathcal Log}^s({\mathcal M}^0_{{\mathcal O}_{\dot{K}}})$ with $G$-action which by Lemma \ref{comalg} is $\dot{\pi}$-adically complete and $\dot{\pi}$-torsion free. Since also ${\mathcal M}^0_{{\mathcal O}_{\dot{K}}}\otimes_{{\mathcal O}_K}\Omega^{s}_{{\mathfrak{X}}}$ is $\dot{\pi}$-adically complete and $\dot{\pi}$-torsion free we have a $G$-equivariant map\begin{gather}{\mathcal Log}^s({\mathcal M}^0_{{\mathcal O}_{\dot{K}}})\longrightarrow {\mathcal M}^0_{{\mathcal O}_{\dot{K}}}\otimes_{{\mathcal O}_K}\Omega^{s}_{{\mathfrak{X}}}.\label{keymap}\end{gather}Recall that we view the $k$-scheme $Y$ from section \ref{cohmodp} as (the central) irreducible component of ${\mathfrak X}\otimes k$; in this way the open subscheme $Y^0\subset Y$ is also open in ${\mathfrak X}\otimes k$.

\begin{lem}\label{lofore} ${\mathcal Log}^s({\mathcal M}^0_{{\mathcal O}_{\dot{K}}})\otimes k$ is a $G$-equivariant subsheaf of $({\mathcal M}^0_{{\mathcal O}_{\dot{K}}}\otimes_{{\mathcal O}_K}\Omega^{s}_{{\mathfrak{X}}})\otimes_{{\mathcal O}_{\dot{K}}}k$ which on $Y^0$ restricts to $\widetilde{{M}}\otimes\underline{\mathbb L}^s_Y|_{Y^0}$.
\end{lem}

{\sc Proof:} From Lemma \ref{comalg} (b) we know that the inclusion ${\mathcal Log}_{alg}^s({\mathcal M}^0_{{\mathcal O}_{\dot{K}}})\to{\mathcal M}^0_{{\mathcal O}_{\dot{K}}}\otimes_{{\mathcal O}_K}\Omega^{s}_{{\mathfrak{X}}}$ induces an injective map of sheaves\begin{gather}{\mathcal Log}_{alg}^s({\mathcal M}^0_{{\mathcal O}_{\dot{K}}})\otimes_{{\mathcal O}_{\dot{K}}}k\longrightarrow({\mathcal M}^0_{{\mathcal O}_{\dot{K}}}\otimes_{{\mathcal O}_K}\Omega^{s}_{{\mathfrak{X}}})\otimes_{{\mathcal O}_{\dot{K}}}k.\label{injre}\end{gather}From Lemma \ref{comalg} (a) we know that the map$${\mathcal Log}_{alg}^s({\mathcal M}^0_{{\mathcal O}_{\dot{K}}})\otimes_{{\mathcal O}_{\dot{K}}}k\longrightarrow{\mathcal Log}^s({\mathcal M}^0_{{\mathcal O}_{\dot{K}}})\otimes_{{\mathcal O}_{\dot{K}}}k$$is an epimorphism of sheaves. Together we conclude that the natural map$${\mathcal Log}^s({\mathcal M}^0_{{\mathcal O}_{\dot{K}}})\otimes_{{\mathcal O}_{\dot{K}}}k\longrightarrow({\mathcal M}^0_{{\mathcal O}_{\dot{K}}}\otimes_{{\mathcal O}_K}\Omega^{s}_{{\mathfrak{X}}})\otimes_{{\mathcal O}_{\dot{K}}}k$$is injective and that its image is the same as that of (\ref{injre}). We now prove our statement concerning the restriction to $Y^0$ of this image sheaf. Since ${\mathcal M}^0_{{\mathcal O}_{\dot{K}}}|_{Y^0}$ is the constant sheaf generated by the free ${\mathcal O}_{\dot{K}}$-module $M^0\otimes_{{\mathcal O}_{{K}}}{\mathcal O}_{\dot{K}}$ with reduction $\widetilde{M}=M^0/\pi.M^0$, it is clear that we may assume $M=K$, the trivial representation. What we must show then is$$({\mathcal Log}^s\bigcap\Gamma(Y^0,\Omega^{s}_{\mathfrak X}))\otimes k={\mathbb{L}}_Y^{s}$$for all $s\ge0$. For $s=0$ both sides are identified with $k$, and the case $s>1$ is reduced to the case $s=1$ by taking exterior products. Thus we assume $s=1$. The containment of the left hand side in the right hand side is clear. Let now $z^{n}\dlog(z)$ be a typical generator of ${\mathbb{L}}_Y^{1}$ as in equation (\ref{frofred}); we need to show that it lies in $({\mathcal Log}^1\bigcap\Gamma(Y^0,\Omega^{1}_{\mathfrak X}))\otimes k$ (here $z=y_1$ for $y_1,\ldots,y_d$ as in equation (\ref{frofred})). The case $n=0$ is clear, and the case $n<0$ is reduced to the case $n>0$ observing $\dlog(z)=-\dlog(z^{-1})$), thus we assume $n>0$. We lift the system $y_1,\ldots,y_d$ on $Y$ to a system $y_1,\ldots,y_d$ on $X$ as in equation (\ref{frofo}), and we also write $z=y_1$ for the lifted $y_1$. Choose pairwise distinct $a_0,\ldots,a_{n}\in{\mathcal{O}}_K$. Since the matrix $(a_i^j)_{0\le i,j\le n}$ is invertible over $K$ (Vandermonde) we may find $x_0,\ldots,x_n\in{\mathcal{O}}_K$ such that, if we set $c_j=\sum_{i=0}^{n}x_ia_i^j$, then $c_j=0$ for $0\le j<n$ and $c_n\ne0$ (possibly a very small $c_n$ since $(a_i^j)_{0\le i,j\le n}$ may not be invertible over ${\mathcal{O}}_K$). Write $c_j=\sum_{i=0}^{n}x_ia_i^j$ for any $j\ge 0$. For $m\in\mathbb{N}$ we have$$\sum_{i=0}^n\frac{x_i}{1-a_i\pi^mz}=\sum_{j=0}^{\infty}c_j\pi^{mj}z^j.$$Now fix $m\in \mathbb{N}$ such that $|\pi^mc_j|<|c_n|$ for all $j>n$ with $|c_j|>|c_n|$. Then $|c_j\pi^{mj}|<|c_n\pi^{mn}|$ for all $j>n$. Hence$$(c_n\pi^{mn})^{-1}\sum_{i=0}^nx_i\dlog(1-a_i\pi^mz)\in{\mathcal Log}^1\bigcap\Gamma(Y^0,\Omega^{1}_{\mathfrak X})$$lifts the form $z^{n}\dlog(z)\in{\mathbb{L}}_Y^{1}$.\hfill$\Box$\\

For $j, t\in\{0,\ldots,d\}$ let$$F^j(t)=\{Z\in F^j\quad|\quad Z=Z_0\cap\ldots\cap Z_j\mbox{ with }Z_i\in F^0(t)\mbox{ for all }0\le i\le j\}.$$Note that $F^j(t)$ is stable under ${\rm SL}_{d+1}(K)$ (because $F^0(t)$ is stable under ${\rm SL}_{d+1}(K)$). For any $t\in\{0,\ldots,d\}$ with $F^0(t)\ne \emptyset$ (i.e. with $t\equiv n|M|$ modulo $(d+1)$ for some $n\in\mathbb{Z}$), the minimal number $j$ with $F^j(t)=\emptyset$ is the quotient of $d+1$ divided by the order of (the class of) $|M|$ in $\mathbb{Z}/(d+1)$ (we set $F^{d+1}(t)=\emptyset$).

\begin{satz}\label{logred} There is a canonical assignment which to $j, t, s\in\{0,\ldots,d\}$ and $Z\in F^j(t)$ assigns a sheaf ${\mathbb L}^{s}(\widetilde{\mathcal M})_Z$ on ${\mathfrak{X}}\otimes k$ with the following properties. ${\mathbb L}^{s}(\widetilde{\mathcal M})_Z$ is supported on $Z$. For $Z=Y\in F^0(0)$ we have ${\mathbb L}^{s}(\widetilde{\mathcal M})_Y=\iota_{Y,*}\widetilde{\mathcal M}[{\mathbb L}^{s}_Y]$ (as defined earlier). There is a ${\rm SL}_{d+1}(K)$-stable direct sum decomposition$${\mathcal Log}^s({\mathcal M}^0_{{\mathcal O}_{\dot{K}}})\otimes k\cong\bigoplus_{t\in\{0,\ldots,d\}}({\mathcal Log}^s({\mathcal M}^0_{{\mathcal O}_{\dot{K}}})\otimes k)(t)$$and for each $t\in\{0,\ldots,d\}$ a ${\rm SL}_{d+1}(K)$-equivariant long exact sequence\begin{gather}0\longrightarrow ({\mathcal Log}^s({\mathcal M}^0_{{\mathcal O}_{\dot{K}}})\otimes k)(t)\longrightarrow\bigoplus_{Z\in F^0(t)}{\mathbb L}^{s}(\widetilde{\mathcal M})_Z\longrightarrow\bigoplus_{Z\in F^1(t)}{\mathbb L}^{s}(\widetilde{\mathcal M})_Z\longrightarrow\ldots.\label{seqlofo}\end{gather}
\end{satz}

{\sc Proof:} The direct sum decomposition of $({\mathcal M}^0_{{\mathcal O}_{\dot{K}}}\otimes_{{\mathcal O}_K}{\mathcal O}_{{\mathfrak{X}}})\otimes_{{\mathcal O}_{\dot{K}}}k$  from Theorem \ref{strumn} yields the analoguous decomposition $$({\mathcal M}^0_{{\mathcal O}_{\dot{K}}}\otimes_{{\mathcal O}_K}\Omega^s_{{\mathfrak{X}}})\otimes_{{\mathcal O}_{\dot{K}}}k\cong\bigoplus_{t\in\{0,\ldots,d\}}\frac{({\mathcal M}^0_{{\mathcal O}_{\dot{K}}}\otimes_{{\mathcal O}_K}\Omega^s_{{\mathfrak{X}}})\otimes_{{\mathcal O}_{\dot{\mathfrak{X}}}}{\mathcal O}_{\overline{\mathfrak{X}(t)}}}{{\mathcal O}_{\overline{\mathfrak{X}(t)}}\mbox{-torsion}}$$where the summand for $t\in\{0,\ldots,d\}$ is locally free on $\overline{\mathfrak{X}(t)}$. We have$$\frac{({\mathcal M}^0_{{\mathcal O}_{\dot{K}}}\otimes_{{\mathcal O}_K}\Omega^s_{{\mathfrak{X}}})\otimes_{{\mathcal O}_{\dot{\mathfrak{X}}}}{\mathcal O}_{\overline{\mathfrak{X}(0)}}}{{\mathcal O}_{\overline{\mathfrak{X}(0)}}\mbox{-torsion}}\otimes_{{\mathcal O}_{\overline{\mathfrak{X}(0)}}}{\mathcal O}_Y=\iota_{Y,*}(\widetilde{\mathcal M}[{\mathcal O}_Y]\otimes_{{\mathcal O}_Y}\Omega^s_Y).$$Now $\iota_{Y,*}\widetilde{\mathcal M}[{\mathbb L}^{s}_Y]\subset \iota_{Y,*}(\widetilde{\mathcal M}[{\mathcal O}_Y]\otimes_{{\mathcal O}_Y}\Omega^s_Y)$ by the definition of $\widetilde{\mathcal M}[{\mathbb L}^{s}_Y]$. We let ${\mathbb L}^{s}(\widetilde{\mathcal M})_Y=\iota_{Y,*}\widetilde{\mathcal M}[{\mathbb L}^{s}_Y]$ and then we move this definition around by means of the action of $G$ to obtain for each $Z\in F^0(t)$ (any $t$) a subsheaf $${\mathbb L}^{s}(\widetilde{\mathcal M})_Z\subset \frac{({\mathcal M}^0_{{\mathcal O}_{\dot{K}}}\otimes_{{\mathcal O}_K}\Omega^s_{{\mathfrak{X}}})\otimes_{{\mathcal O}_{\dot{\mathfrak{X}}}}{\mathcal O}_{\overline{\mathfrak{X}(t)}}}{{\mathcal O}_{\overline{\mathfrak{X}(t)}}\mbox{-torsion}}\otimes_{{\mathcal O}_{\overline{\mathfrak{X}(t)}}}{\mathcal O}_Z.$$ We have ${\mathbb L}^{s}(\widetilde{\mathcal M})_Z|_{Z\cap Z'}={\mathbb L}^{s}(\widetilde{\mathcal M})_{Z'}|_{Z\cap Z'}$ for all $Z,Z'\in F^0(t)$ (because of $G$-equivariance: there are $g\in G$ which interchange $Z$ and $Z'$). This means that also for $j>0$ we obtain subsheaves$${\mathbb L}^{s}(\widetilde{\mathcal M})_Z\subset \frac{({\mathcal M}^0_{{\mathcal O}_{\dot{K}}}\otimes_{{\mathcal O}_K}\Omega^s_{{\mathfrak{X}}})\otimes_{{\mathcal O}_{\dot{\mathfrak{X}}}}{\mathcal O}_{\overline{\mathfrak{X}(t)}}}{{\mathcal O}_{\overline{\mathfrak{X}(t)}}\mbox{-torsion}}\otimes_{{\mathcal O}_{\overline{\mathfrak{X}(t)}}}{\mathcal O}_Z$$for each $Z\in F^j(t)$ and ${\rm SL}_{d+1}(K)$-stable subsheaves $${\mathcal F}(t)\subset \frac{({\mathcal M}^0_{{\mathcal O}_{\dot{K}}}\otimes_{{\mathcal O}_K}\Omega^s_{{\mathfrak{X}}})\otimes_{{\mathcal O}_{\dot{\mathfrak{X}}}}{\mathcal O}_{\overline{\mathfrak{X}(t)}}}{{\mathcal O}_{\overline{\mathfrak{X}(t)}}\mbox{-torsion}}$$such that there are long exact sequences$$0\longrightarrow {\mathcal F}(t)\longrightarrow\bigoplus_{Z\in F^0(t)}{\mathbb L}^{s}(\widetilde{\mathcal M})_Z\longrightarrow\bigoplus_{Z\in F^1(t)}{\mathbb L}^{s}(\widetilde{\mathcal M})_Z\longrightarrow\ldots.$$The restriction of ${\mathbb L}^{s}(\widetilde{\mathcal M})_Y=\iota_{Y,*}\widetilde{\mathcal M}[{\mathbb L}^{s}_Y]$ and hence of ${\mathcal F}(0)$ to the open subscheme $Y^0$ is just $\widetilde{{M}}\otimes\underline{\mathbb L}^s_Y|_{Y^0}$. In view of Lemma \ref{lofore} and $G$-equivariance we conclude that the subsheaves ${\mathcal Log}^s({\mathcal M}^0_{{\mathcal O}_{\dot{K}}})\otimes k$ and $\oplus_{t\in\{0,\ldots,d\}}{\mathcal F}(t)$ of $({\mathcal M}^0_{{\mathcal O}_{\dot{K}}}\otimes_{{\mathcal O}_K}\Omega^s_{{\mathfrak{X}}})\otimes_{{\mathcal O}_{\dot{K}}}k$ coincide when restricted to $Y^0$ and to each $G$-translate of $Y^0$ in ${\mathfrak{X}}\otimes k$. By their construction both these subsheaves are maximal inside $({\mathcal M}^0_{{\mathcal O}_{\dot{K}}}\otimes_{{\mathcal O}_K}\Omega^s_{{\mathfrak{X}}})\otimes_{{\mathcal O}_{\dot{K}}}k$ with this given restriction to $Y^0$ and its $G$-translates, hence they coincide.\hfill$\Box$\\

\begin{pro}\label{globlo} If the weights of $M$ are small then the map (\ref{keymap}) induces for any $i$ an isomorphism$$H^i({\mathfrak{X}},{\mathcal Log}^s({\mathcal M}^0_{{\mathcal O}_{\dot{K}}}))\cong H^i({\mathfrak{X}},{\mathcal M}^0_{{\mathcal O}_{\dot{K}}}\otimes_{{\mathcal O}_K}\Omega^{s}_{{\mathfrak{X}}}).$$
\end{pro}

{\sc Proof:} For the sheaves ${\mathcal F}={\mathcal Log}^s({\mathcal M}^0_{{\mathcal O}_{\dot{K}}})$ and ${\mathcal F}={\mathcal M}^0_{{\mathcal O}_{\dot{K}}}\otimes_{{\mathcal O}_K}\Omega^{s}_{{\mathfrak{X}}}$ we have the spectral sequences\begin{gather} E_{2}^{pq}=R^p\lim_{\leftarrow\atop m}(H^q({\mathfrak X},{\mathcal F}_m))\Rightarrow H^{p+q}({\mathfrak X},\lim_{\leftarrow\atop m}{\mathcal F}_m)=H^{p+q}({\mathfrak X},{\mathcal F})\notag\end{gather}where $(.)_m$ denotes reduction modulo $\dot{\pi}^m$. The map (\ref{keymap}) induces a map between these spectral sequences and we see that it is enough to show\begin{gather}H^i({\mathfrak{X}},({\mathcal Log}^s({\mathcal M}^0_{{\mathcal O}_{\dot{K}}}))_m)\cong H^i({\mathfrak{X}},({\mathcal M}^0_{{\mathcal O}_{\dot{K}}}\otimes_{{\mathcal O}_K}\Omega^{s}_{{\mathfrak{X}}})_m)\label{modem}\end{gather}for any $m\ge 1$, any $i\ge0$. Since ${\mathcal F}={\mathcal Log}^s({\mathcal M}^0_{{\mathcal O}_{\dot{K}}})$ and ${\mathcal F}={\mathcal M}^0_{{\mathcal O}_{\dot{K}}}\otimes_{{\mathcal O}_K}\Omega^{s}_{{\mathfrak{X}}}$ are ${\mathcal O}_{\dot{K}}$-flat we get exact sequences of sheaves$$0\to{\mathcal F}_{m-1}\stackrel{\dot{\pi}^{m-1}}{\longrightarrow}{\mathcal F}_m\longrightarrow{\mathcal F}_1\longrightarrow0.$$Comparing the associated long exact cohomology sequences we reduce our task to proving the isomorphism (\ref{modem}) in the case $m=1$, i.e. to proving$$H^*({\mathfrak{X}},{\mathcal Log}^s({\mathcal M}^0_{{\mathcal O}_{\dot{K}}})\otimes k)\cong H^*({\mathfrak{X}},{\mathcal M}^0_{{\mathcal O}_{\dot{K}}}\otimes_{{\mathcal O}_K}\Omega^{s}_{{\mathfrak{X}}}\otimes k).$$First suppose $|M|\notin(d+1).\mathbb{Z}$. Then our hypothesis that the weights of $M$ be small implies that the order of (the class of) $|M|$ in $\mathbb{Z}/(d+1)$ is $d+1$, cf. the proof of Lemma \ref{wesm}. Then comparing Theorem \ref{logred} with the result from Theorem \ref{strumn} and using $G$-equivariance we reduce to proving$$H^i(Y,\widetilde{\mathcal M}[{\mathbb L}^{s}_Y])\cong H^i(Y,{\widetilde{\mathcal M}}[{\mathcal O}_Y]\otimes_{{\mathcal O}_Y}\Omega^s_Y)$$for any $i$. But this we did in Theorem \ref{klewe}. Now suppose $|M|\in(d+1).\mathbb{Z}$. Under our hypothesis that the weights of $M$ be small this means $M|_{{\rm SL}_{d+1}(K)}$ is trivial, hence ${\mathcal M}^0_{{\mathcal O}_{\dot{K}}}$ is the constant sheaf with value ${\mathcal O}_{\dot{K}}$. Since $\Omega^{s}_{{\mathfrak{X}}}\otimes k$ is locally free over ${\mathcal O}_{{\mathfrak{X}}\otimes k}$ we have an exact sequence$$0\longrightarrow\Omega^{s}_{{\mathfrak{X}}}\otimes k\longrightarrow\bigoplus_{Z\in F^0}\Omega^{s}_{{\mathfrak{X}}}\otimes{\mathcal O}_Z\longrightarrow\bigoplus_{Z\in F^1}\Omega^{s}_{{\mathfrak{X}}}\otimes{\mathcal O}_Z\longrightarrow\ldots.$$On the other hand the exact sequence (\ref{seqlofo}) becomes $$0\longrightarrow {\mathcal Log}^s(\underline{{\mathcal O}_{\dot{K}}})\otimes k\longrightarrow\bigoplus_{Z\in F^0}{\mathbb L}^{s}(\underline{k})_Z\longrightarrow\bigoplus_{Z\in F^1}{\mathbb L}^{s}(\underline{k})_Z\longrightarrow\ldots$$with each ${\mathbb L}^{s}(\underline{k})_Z$ the push forward to ${\mathfrak X}$ of a constant sheaf on $Z$ (which we denote by ${\mathbb L}^{s}(\underline{k})_Z$, too). Comparing we reduce to proving$$H^*(Z,{\mathbb L}^{s}(\underline{k})_Z)\cong H^*(Z,\Omega^s_{{\mathfrak{X}}}\otimes{\mathcal O}_Z)$$for any $Z\in F^j$, any $j$. By $G$-equivariance we may ssume $Z\subset Y$. Let $\iota_Y^Z:Z\to Y$ denote the closed embedding. The proof of Theorem \ref{logred} shows $(\iota_Y^Z)_*{\mathbb L}^{s}(\underline{k})_Z={\mathbb L}^{s}_Z(0)$ as defined in Theorem \ref{konskomop} (b). Hence we may conclude by that Theorem.\hfill$\Box$\\

{\it Remarks:} (1) From Proposition \ref{globlo} it follows (take $M=K$) that every bounded differential $s$-form on $X$, i.e. every element of $H^0({\mathfrak{X}},\Omega^{s}_{{\mathfrak{X}}})\otimes K$, is in fact logarithmic, in particular it is closed. Thus $H^0({\mathfrak{X}},\Omega^{s}_{{\mathfrak{X}}})\otimes K$ must be the space of bounded logarithmic differential $s$-forms on $X$ studied in \cite{iovspi} (if $\kara(K)=0$).\\
(2) Suppose $M=K$, the trivial $G$-representation. Let $W\omega^{\bullet}_{\mathfrak{X}}$ denote the logarithmic de Rham complex of the special fibre of ${\mathfrak{X}}$. The same proof as for \ref{globlo} provides isomorphisms $$H^j({\mathfrak{X}},{\mathcal Log}^{\bullet}({\mathcal M}^0_{{\mathcal O}_{\dot{K}}}))\cong H^j({\mathfrak{X}},W\omega^{\bullet}_{\mathfrak{X}})$$for any $j$, hence altogether isomorphisms$$H^j({\mathfrak{X}},\Omega^{{\bullet}}_{{\mathfrak{X}}})\cong H^j({\mathfrak{X}},W\omega^{\bullet}_{\mathfrak{X}}).$$Similarly, for the quotients ${\mathfrak{X}}_{\Gamma}$ of ${\mathfrak{X}}$ as in section \ref{norehose} we get by the same proof isomorphisms\begin{gather}H^j({\mathfrak{X}}_{\Gamma},\Omega^{{\bullet}}_{{\mathfrak{X}_{\Gamma}}})\cong H^j({\mathfrak{X}}_{\Gamma},W\omega^{\bullet}_{\mathfrak{X}_{\Gamma}}).\label{luc}\end{gather}These isomorphisms (\ref{luc}) are those constructed by Hyodo (see \cite{illast}) by means of $p$-adic \'{e}tale sheaves of vanishing cycles for general projective semistable schemes with ordinary reduction. They must not be confused with the Hyodo-Kato isomorphisms which are used to define the filtered $(\phi,N)$-modules which recover the $p$-adic \'{e}tale cohomology of the generic fibre of ${\mathfrak{X}}_{\Gamma}$. 

\section{The Hodge spectral sequence}

\label{norehose}

Let $\Gamma\subset {\rm SL}\sb {d+1}(K)$ be a discrete torsionfree and cocompact subgroup. It is proved in \cite{mus} that the quotient ${\mathfrak{X}}_{\Gamma}=\Gamma\backslash{\mathfrak{X}}$ is the $\pi$-adic formal completion of a projective ${\mathcal O}_K$-scheme. Passing to a smaller $\Gamma$ if necessary we may assume that ${\mathfrak{X}}_{\Gamma}$ has strictly semistable reduction, i.e. all irreducible components of ${\mathfrak{X}}_{\Gamma}$ are smooth. Let $X_{\Gamma}=\Gamma\backslash X={\mathfrak{X}}_{\Gamma}\otimes K$, the analytification of a smooth projective $K$-scheme. Let $M$ be a $K[\Gamma]$-module with $\dim_KM<\infty$; we write $\underline{M}={\mathcal M}$ for the constant sheaf on $X$, resp. on $\mathfrak{X}$, generated by $M$. 

For a $\Gamma$-equivariant sheaf ${\mathcal F}$ on ${\mathfrak{X}}$ or $X$ we write ${\mathcal F}^{\Gamma}$ for the descended sheaf on ${\mathfrak{X}}_{\Gamma}$ or $X_{\Gamma}$. For example, the constant local system $\underline{M}={\mathcal M}$ on $\mathfrak{X}$ or $X$ gives rise to a (non constant in general !) descended local system ${\mathcal {M}}^{\Gamma}$ on $\mathfrak{X}_{\Gamma}$ or $X_{\Gamma}$. We are interested in the cohomology of the sheaf complex ${\mathcal {M}}^{\Gamma}\otimes_K\Omega^{\bullet}_{{X}_{\Gamma}}=(\underline{M}\otimes_K\Omega^{\bullet}_{X})^{\Gamma}$ on $\mathfrak{X}_{\Gamma}$ or $X_{\Gamma}$. The Hodge spectral sequence\begin{gather}E_1^{r,s}=H^s(X_{\Gamma},{\mathcal {M}}^{\Gamma}\otimes_K\Omega^{r}_{X_{\Gamma}} )\Rightarrow H^{r+s}(X_{\Gamma},{\mathcal {M}}^{\Gamma}\otimes_K\Omega^{\bullet}_{{X}_{\Gamma}}    )
\label{hodss}\end{gather}gives rise to the Hodge filtration$$H^{d}(X_{\Gamma},{\mathcal {M}}^{\Gamma}\otimes_K\Omega^{\bullet}_{{X}_{\Gamma}}    )=F^0_{H}\supset F^1_{H}\supset\ldots\supset F^{d+1}_{H}=0.$$

\begin{kor}\label{hosssp} If $M$ is a $K$-rational $G$-representation with small weights then the Hodge spectral sequence (\ref{hodss}) degenerates in $E_1$. The Hodge filtration $F_H^{\bullet}$ has a canonical splitting defined through logarithmic differential forms (at least after base extension $K\to\dot{K}$).
\end{kor}

{\sc Proof:} We may extend scalars $K\to\dot{K}$. We continue to use the same names for coherent sheaves on $X_{\Gamma}$ and for their push forward to ${\mathfrak X}_{\Gamma}$. We have an inclusion of sheaf complexes ${\mathcal Log}^{\bullet}({\mathcal M}^0_{{\mathcal O}_{\dot{K}}})^{\Gamma}\otimes_{{\mathcal O}_{\dot{K}}} \dot{K}\to{\mathcal {M}}^{\Gamma}\otimes_K\Omega^{\bullet}_{{X}_{\Gamma}}\otimes_K\dot{K}$ on ${\mathfrak X}_{\Gamma}$ with trivial differentials on the former. Therefore it is enough to prove that for any $0\le s\le d$ the natural maps$$H^*({\mathfrak X}_{\Gamma},{\mathcal Log}^{s}({\mathcal M}^0_{{\mathcal O}_{\dot{K}}})^{\Gamma}\otimes_{{\mathcal O}_{\dot{K}}} \dot{K})\longrightarrow H^*({\mathfrak X}_{\Gamma},{\mathcal {M}}^{\Gamma}\otimes_K\Omega^{s}_{{X}_{\Gamma}}\otimes_K\dot{K})$$are isomorphisms.  Now $${\mathcal {M}}^{\Gamma}\otimes_K\Omega^{s}_{{X}_{\Gamma}}\otimes_K{\dot{K}}=({\mathcal M}^0_{{\mathcal O}_{\dot{K}}}\otimes_{{\mathcal O}_K}\Omega^{s}_{{\mathfrak{X}}})^{\Gamma}\otimes_{{\mathcal O}_{\dot{K}}}{\dot{K}}.$$Since ${\mathfrak X}_{\Gamma}$ is quasicompact, taking cohomology commutes with applying $(.)\otimes_{{\mathcal O}_{\dot{K}}}{\dot{K}}$. Therefore it will be enough to show$$H^*({\mathfrak X}_{\Gamma},{\mathcal Log}^s({\mathcal M}^0_{{\mathcal O}_{\dot{K}}})^{\Gamma})\cong H^*({\mathfrak X}_{\Gamma},({\mathcal M}^0_{{\mathcal O}_{\dot{K}}}\otimes_{{\mathcal O}_K}\Omega^{s}_{{\mathfrak{X}}})^{\Gamma}).$$For both ${\mathcal F}={\mathcal Log}^s({\mathcal M}^0_{{\mathcal O}_{\dot{K}}})$ and ${\mathcal F}={\mathcal M}^0_{{\mathcal O}_{\dot{K}}}\otimes_{{\mathcal O}_K}\Omega^{s}_{{\mathfrak{X}}}$ we have the spectral sequence$$E_2^{rt}=H^r(\Gamma, H^t({\mathfrak{X}},{\mathcal F}))\Rightarrow H^{r+t}({\mathfrak{X}}_{\Gamma},{\mathcal F}^{\Gamma}).$$
We conclude by Proposition \ref{globlo} (alternatively we could repeat the proof of Proposition \ref{globlo}). \hfill$\Box$\\ 

Let again $M$ be an arbitrary $K[\Gamma]$-module with $\dim_KM<\infty$. From now on we suppose $\kara(K)=0$. For an open subscheme $U$ of $\mathfrak{X}\otimes k$ we denote by $\overline{U}$ the Zariski closure of $U$ in ${\mathfrak{X}}\otimes k$, and by $]\overline{U}[=]\overline{U}[_{\mathfrak{X}}=sp^{-1}(\overline{U})$ its tube in $X$, the preimage under the specialization map $sp:X\to\mathfrak{X}\otimes k$. For $i\ge0$ we define the sheaf ${\mathbb{L}}^i(M)$ on ${\mathfrak{X}}\otimes k$ (or equivalently: on ${\mathfrak{X}}$) by setting\footnote{Logically the notation ${\mathbb{L}}$ as used here has nothing to do with the notation ${\mathbb{L}}$ as used in the previous sections; however, the ${\mathbb{L}}$'s play the same role in their respective contexts.}$${\mathbb{L}}^i(M)(U)=\ke[\underline{M}\otimes\Omega^{i}_{X}(]\overline{U}[)\longrightarrow\underline{M}\otimes\Omega^{i+1}_{X}(]\overline{U}[)]$$for open $U\subset{\mathfrak{X}}\otimes k$. We get a sheaf complex ${\mathbb{L}}^{\bullet}(M)$ on ${\mathfrak{X}}$ with trivial differentials. For $i\ge0$ let $\tau_i(\underline{M}\otimes\Omega^{\bullet}_X)$ be the subsheaf complex of $\underline{M}\otimes\Omega^{\bullet}_{X}$ on ${\mathfrak{X}}\otimes k$ whose value $\tau_i(\underline{M}\otimes\Omega^{\bullet}_X)(U)$  for open $U\subset{\mathfrak{X}}\otimes k$ is the complex $$\underline{M}\otimes\Omega^0_X(]\overline{U}[)\longrightarrow\ldots\longrightarrow\underline{M}\otimes\Omega^{i-1}_X(]\overline{U}[)\longrightarrow {\mathbb{L}}^i(M)(U)\longrightarrow0\longrightarrow\ldots.$$

We write $\tau_i({{\mathcal {M}}}^{\Gamma}\otimes_K\Omega^{\bullet}_{{X}_{\Gamma}})=\tau_i({\underline{M}}\otimes_K\Omega^{\bullet}_{X})^{\Gamma}$ for the descended sheaf complex on $\mathfrak{X}_{\Gamma}$ or $X_{\Gamma}$.

 For a complex $C^{\bullet}=(C^0\stackrel{d}{\to}C^1\stackrel{d}{\to}C^2\stackrel{d}{\to} \ldots)$ (of abstract groups, or sheaves) we put$$t_{\le i}C^{\bullet}=(C^0\stackrel{d}{\longrightarrow}\ldots\stackrel{d}{\longrightarrow}C^{i-1}\stackrel{d}{\longrightarrow}\ke(d)\stackrel{d}{\longrightarrow}0\stackrel{d}{\longrightarrow}\ldots).$$

\begin{pro}\label{acycanw} We have$$H^t({\mathfrak{X}},\tau_i({\underline{M}}\otimes_K\Omega^{\bullet}_{{X}}))= \left\{\begin{array}{l@{\quad:\quad}l}M\otimes_K H_{dR}^t(X)&0\le t\le i\\0&t>i\end{array}\right..$$In particular,$$H^d({\mathfrak{X}}_{\Gamma},\tau_i({{\mathcal {M}}}^{\Gamma}\otimes_K\Omega^{\bullet}_{{X}_{\Gamma}}))=H^d(\Gamma,M\otimes_Kt_{\le i}\Omega_X^{\bullet}(X)).$$
\end{pro}

{\sc Proof:}We first deduce the second statement from the first one. Since $X$ is a Stein space we have $H^s(X,\Omega_X^r)=0$ for all $r\ge0$, all $s>0$ (see \cite{kiaub}), hence $$H^t_{dR}(X)=H^t(X,\Omega_X^{\bullet})=\frac{t_{\le t}\Omega_X^{\bullet}(X)}{t_{\le t-1}\Omega_X^{\bullet}(X)}[-t]$$for all $t$ (the last term is a complex concentrated in degree $0$). Together with the first statement we deduce that the natural map of sheaf complexes $\tau_i({\underline{M}}\otimes_K\Omega^{\bullet}_{{X}})\to M\otimes_K t_{\le i}\Omega_X^{\bullet}$ induces an isomorphism$${\mathbb R}\Gamma({\mathfrak{X}},\tau_i({\underline{M}}\otimes_K\Omega^{\bullet}_{{X}}))\cong M\otimes_Kt_{\le i}\Omega_X^{\bullet}(X).$$This gives the second statement. The first one will be deduced from de Shalit's acyclicity theorem. We may of course assume $M=K$. It will be enough to show$$H^t({\mathfrak{X}},\frac{\tau_i\Omega^{\bullet}_{{X}}}{\tau_{i-1}\Omega^{\bullet}_{{X}}})=\left\{\begin{array}{l@{\quad:\quad}l}H_{dR}^i(X)&t=i\\0&t\ne i\end{array}\right.$$where we set $\tau_{-1}(\underline{M}\otimes\Omega^{\bullet}_X)=0$. For $T\in F^s$ (any $s$) let $${\mathcal Z}(T)=\{Z\in F^0\quad|\quad T\subset Z\}$$and $\dot{T}=\cup_{Z\in{\mathcal Z}(T)}Z$. For all sufficiently small open neighbourhoods $U\subset{\mathfrak{X}}\otimes k$ of a given closed point of ${\mathfrak{X}}\otimes k$ we have $\overline{U}=\dot{T}$ for some $T$. Then $]\overline{U}[=]\dot{T}[$ is a Stein space, hence$$H^i(]\overline{U}[,\frac{\tau_i\Omega^{\bullet}_{{X}}}{\tau_{i-1}\Omega^{\bullet}_{{X}}})=H_{dR}^i(]\overline{U}[).$$Therefore, if ${\mathcal{H}}_{dR}^i$ denotes the sheaf associated with the presheaf$$U\mapsto H_{dR}^i(]\overline{U}[)$$on ${{\mathfrak{X}}}$, then we must show$$H^t({{\mathfrak{X}}},{\mathcal{H}}_{dR}^i)=\left\{\begin{array}{l@{\quad:\quad}l}H^i_{dR}(X)&t=0\\0&t\ne0\end{array}\right..$$For $T\in F^s$ (any $s$) let $\dot{T}^0$ denote the maximal open subscheme of ${\mathfrak{X}}\otimes k$ which is contained in $\dot{T}$. We compute $H^t({{\mathfrak{X}}},{\mathcal{H}}_{dR}^i)$ as Cech cohomology with respect to the open covering ${{\mathfrak{X}}}=\bigcup_{T\in F^d}\dot{T}^0$. Note that for any collection $(T_1,\ldots,T_{r+1})\in (F^d)^{r+1}$ the intersection $\dot{T}_1^0\cap\ldots\cap\dot{T}^0_{r+1}$ is empty or equals $\dot{T}^0$ for some $T\in F^s$, some $s$. In the latter case it follows that $\overline{\dot{T}_1^0\cap\ldots\cap\dot{T}^0_{r+1}}=\dot{T}$. From the definition of ${\mathcal{H}}_{dR}^i$ we know on the other hand that for all $T\in F^s$, all $s$, we have $H^r(]\dot{T}[,{\mathcal{H}}_{dR}^i)=0$ for all $r>0$. Together we get$$H^r(\dot{T}_1^0\cap\ldots\cap\dot{T}^0_{r+1},{\mathcal{H}}_{dR}^i)=0$$for all $r>0$. Therefore it will be enough to show that the complex$$\prod_{T\in F^d}H^i_{dR}(]\dot{T}[)\longrightarrow\prod_{(T_1,T_2)\in(F^{d})^2}H^i_{dR}(]\dot{T_1}\bigcap\dot{T_2}[)\longrightarrow\ldots$$is a resolution of $H_{dR}^i(X)$. By de Shalit's acyclicity theorem \cite{ds} (see also \cite{acy}) we know that the complex $$\prod_{T\in F^0}H^i_{dR}(]T[)\longrightarrow\prod_{T\in F^1}H^i_{dR}(]T[)\longrightarrow\ldots$$is a resolution of $H_{dR}^i(X)$. Both these complexes map to the total complex of the double complex$$K^{rs}=\prod_{(T_1,\ldots,T_{r+1}),T'}H_{dR}^i(]T'[)$$where the product is taken over all $(T_1,\ldots,T_{r+1})\in(F^{d})^{r+1}$ and all $T'\in F^s$ such that $T_1\cap\ldots\cap T_{r+1}\subset T'$. It will be enough to show that these maps are quasiisomorphisms. It is clear that for fixed $s$ the complex $K^{\bullet s}$ is a resolution of $\prod_{T\in F^s}H^i_{dR}(]T[)$. On the other hand it follows from Lemma \ref{lokazyk} below that for fixed $r$ the complex $K^{r\bullet}$ is a resolution of $$\prod_{(T_1,\ldots,T_{r+1})\in (F^{d})^{r+1}}H^i_{dR}(]\dot{T_1}\bigcap\ldots\bigcap\dot{T}_{r+1}[)$$and this completes the proof of \ref{acycanw}.\hfill$\Box$
 
\begin{lem}\label{lokazyk} (see \cite{acy} Corollary 2.9 (1)) For any $T\in F^s$ (any $s$) the sequence$$0\longrightarrow H_{dR}^i(]\dot{T}[_{\mathfrak{X}})\longrightarrow \prod_{Z\in{\mathcal Z}(T)}H_{dR}^i(]Z[_{\mathfrak{X}})\longrightarrow \prod_{R\subset {\mathcal Z}(T)\atop |R|=2}H_{dR}^i(]\bigcap_{Z\in R}Z[_{\mathfrak{X}})\longrightarrow\ldots$$is exact.\hfill$\Box$\\
\end{lem}

We have the covering spectral sequence\begin{gather}E_2^{r,s}=H^r(\Gamma,M\otimes_KH^s_{dR}(X))\Rightarrow H^{r+s}(X_{\Gamma},{\mathcal {M}}^{\Gamma}\otimes_K\Omega^{\bullet}_{{X}_{\Gamma}}    )\label{covss}\end{gather}which degenerates in $E_2$, as is shown in \cite{schn}. Denote by $$H^{d}(X_{\Gamma},{\mathcal {M}}^{\Gamma}\otimes_K\Omega^{\bullet}_{{X}_{\Gamma}}    )=F^0_{\Gamma}\supset F^1_{\Gamma}\supset\ldots\supset F^{d+1}_{\Gamma}=0$$the filtration on $H^{d}(X_{\Gamma},{\mathcal {M}}^{\Gamma}\otimes_K\Omega^{\bullet}_{{X}_{\Gamma}}    )$ induced by (\ref{covss}) (it turns out that the cohomology in other degrees is not interesting). By \cite{schn} Theorem 2 and Proposition 2, section 1, we have for $i=0,\ldots,d+1$:\begin{gather}\dim_KF^i_{\Gamma}=\left\{\begin{array}{l@{\quad:\quad}l}(d+1-i)\mu(\Gamma,M)&\mbox{}\,\,d\,\,\mbox{is odd or}\,\,2i>d\\(d+1-i)\mu(\Gamma,M)+\dim_KM^{\Gamma}&\mbox{}\,\,d\,\,\mbox{is even and}\,\,2i\le d\end{array}\right.\label{pecomp}\\\mu(\Gamma,M)=\mu(\Gamma,M^*)\label{dimdu}\end{gather}Here $\mu(\Gamma,M)=\dim_KH^d(\Gamma,M)$ and $M^*=\Hom_K(M,K)$ and we must assume $d\ge2$.\\

{\bf Conjecture:} (Schneider) \begin{gather}H^{d}(X_{\Gamma},{\mathcal {M}}^{\Gamma}\otimes_K\Omega^{\bullet}_{{X}_{\Gamma}}    )=F_H^{i+1}\bigoplus F_{\Gamma}^{d-i}\label{notag}\end{gather}for all $0\le i\le d-1$.

\begin{satz}\label{genfro} The following (i) and (ii) are equivalent:\\(i) The map\begin{gather}H^d({\mathfrak X}_{\Gamma},\tau_i({\mathcal {M}}^{\Gamma}\otimes_K\Omega^{\bullet}_{{X}_{\Gamma}}) \bigoplus{\mathcal {M}}^{\Gamma}\otimes_K\Omega_{{{X}}_{\Gamma},\ge i+1}^{\bullet})\to H^d({\mathfrak X}_{\Gamma},{\mathcal {M}}^{\Gamma}\otimes_K\Omega^{\bullet}_{{X}_{\Gamma}} )\label{surzi}\end{gather}and the analogous map for $M^*$ and $d-i$ (instead of $M$ and $i+1$) are surjective.\\(ii) We have the decomposition\begin{gather}H^{d}(X_{\Gamma},{\mathcal {M}}^{\Gamma}\otimes_K\Omega^{\bullet}_{{X}_{\Gamma}}    )=F_H^{i+1}\bigoplus F_{\Gamma}^{d-i}\label{spliteins}\end{gather}and the analogous decomposition for $M^*$ and $d-i$ (instead of $M$ and $i+1$).
\end{satz}

{\sc Proof:} By definition we have
\begin{align}F_H^{i+1}&=\bi[H^{d}(X_{\Gamma},{\mathcal {M}}^{\Gamma}\otimes_K\Omega_{X_{\Gamma},\ge i+1}^{\bullet})\to H^{d}(X_{\Gamma},{\mathcal {M}}^{\Gamma}\otimes_K\Omega^{\bullet}_{{X}_{\Gamma}}    )]\notag \\F_{\Gamma}^{d-i}&=\bi[H^d(\Gamma,M\otimes_Kt_{\le i}\Omega_X^{\bullet}(X))\to H^d(\Gamma,M\otimes_K\Omega_X^{\bullet}(X))]\notag \\{} &=\bi[H^d(\Gamma,M\otimes_Kt_{\le i}\Omega_X^{\bullet}(X))\to H^{d}(X_{\Gamma},{\mathcal {M}}^{\Gamma}\otimes_K\Omega^{\bullet}_{{X}_{\Gamma}}    )]\notag\end{align}(the last equality holds since $X$ is a Stein space). From \ref{acycanw} it then follows that\begin{gather}F_{\Gamma}^{d-i}=\bi[H^d({\mathfrak{X}}_{\Gamma},\tau_i({\mathcal {M}}^{\Gamma}\otimes_K\Omega^{\bullet}_{{X}_{\Gamma}}))\to H^{d}(X_{\Gamma},{\mathcal {M}}^{\Gamma}\otimes_K\Omega^{\bullet}_{{X}_{\Gamma}} )].\label{charga}\end{gather}This shows that (ii) implies (i). Conversely, if (\ref{surzi}) is surjective then (\ref{charga}) shows\begin{gather}H^{d}(X_{\Gamma},{\mathcal {M}}^{\Gamma}\otimes_K\Omega^{\bullet}_{{X}_{\Gamma}}    )=F_H^{i+1}+ F_{\Gamma}^{d-i}\label{sumhalb}\end{gather}and similarly, if the analog of (\ref{surzi}) with $M^*$ and $d-i$ (instead of $M$ and $i+1$) is surjective then the analog of (\ref{charga}) with $M^*$ and $d-i$ (instead of $M$ and $i+1$) shows the analog of (\ref{sumhalb}) with $M^*$ and $d-i$ (instead of $M$ and $i+1$). By a formal duality argument one then concludes that the sum in (\ref{sumhalb}) is in fact direct. This argument is easily extracted from the proof of \cite{iovspi} Theorem 5.4 and is worked out in a completely analogous situation in the proof of \ref{intelde} below. It rests on Serre duality on the smooth projective $K$-scheme underlying $X_{\Gamma}$ and the computations (\ref{pecomp}) and (\ref{dimdu}) of $\dim_KF^j_{\Gamma}$.\hfill$\Box$\\

{\bf Remark:} As we just saw, the surjectivity of (\ref{surzi}) alone implies (\ref{charga}). This is the sheaf cohomology analog of \cite{schn} p.631, Lemma 2 (ii). To ask in addition for the surjectivity of the analog of (\ref{surzi}) for $M^*$ and $d-i$ for obtaining $F_H^{i+1}\cap F_{\Gamma}^{d-i}=0$ is the strategy of \cite{iovspi}, an alternative to the strategy \cite{schn} p.631, Lemma 2 (i). \\

The inclusion of sheaf complexes ${\mathbb{L}}^{\bullet}(M)^{\Gamma}\to{\mathcal {M}}^{\Gamma}\otimes_K\Omega^{\bullet}_{{X}_{\Gamma}} $ induces a map$$\nabla(M):H^d({\mathfrak X}_{\Gamma},{\mathbb{L}}^{\bullet}(M)^{\Gamma})\longrightarrow H^d({\mathfrak X}_{\Gamma},{\mathcal {M}}^{\Gamma}\otimes_K\Omega^{\bullet}_{{X}_{\Gamma}} ).$$

\begin{kor}\label{nabkri} If $\nabla(M)$ and $\nabla(M^*)$ are surjective then (\ref{spliteins}) holds for all $0\le i\le d-1$.
\end{kor}

{\sc Proof:} The differential in the complex ${\mathbb{L}}^{\bullet}(M)^{\Gamma}$ is zero, consequently the inclusion$${\mathbb{L}}^{\bullet}(M)^{\Gamma}\longrightarrow\tau_i({\mathcal {M}}^{\Gamma}\otimes_K\Omega^{\bullet}_{{X}_{\Gamma}}) \bigoplus{\mathcal {M}}^{\Gamma}\otimes_K\Omega_{{{X}}_{\Gamma},\ge i+1}^{\bullet}$$is a morphism of complexes and \ref{genfro} proves the corollary.\hfill$\Box$\\

By \cite{mus} we know that ${X}_{\Gamma}$ is the analytification of a projective $K$-scheme ${X}_{\Gamma,alg}$. Similarly it follows from GAGA-theorems that the de Rham complex ${\mathcal {M}}^{\Gamma}\otimes_K\Omega^{\bullet}_{{X}_{\Gamma}}$ on ${X}_{\Gamma}$ is the analytification of a complex $({\mathcal {M}}^{\Gamma}\otimes_K\Omega^{\bullet}_{{X}_{\Gamma}})^{alg}$ on ${X}_{\Gamma,alg}$. Consider the conjugate spectral sequence$$E_2^{pq}=H^p({X}_{\Gamma,alg},{\mathcal H}^q(({\mathcal {M}}^{\Gamma}\otimes_K\Omega^{\bullet}_{{X}_{\Gamma}})^{alg}))$$$$\Longrightarrow H^{p+q}({X}_{\Gamma,alg},({\mathcal {M}}^{\Gamma}\otimes_K\Omega^{\bullet}_{{X}_{\Gamma}})^{alg})=H^{p+q}({X}_{\Gamma},{\mathcal {M}}^{\Gamma}\otimes_K\Omega^{\bullet}_{{X}_{\Gamma}}).$$It gives rise to the conjugate filtration$$H^{d}(X_{\Gamma},{\mathcal {M}}^{\Gamma}\otimes_K\Omega^{\bullet}_{{X}_{\Gamma}}    )=H^d({X}_{\Gamma,alg},({\mathcal {M}}^{\Gamma}\otimes_K\Omega^{\bullet}_{{X}_{\Gamma}})^{alg})=F^0_{con}\supset F^1_{con}\supset\ldots\supset F^{d+1}_{con}=0.$$

\begin{pro}\label{conj} Assume\begin{gather}H^{d}(X_{\Gamma},{\mathcal {M}}^{\Gamma}\otimes_K\Omega^{\bullet}_{{X}_{\Gamma}}    )=F_H^{i+1}+F_{con}^{d-i}\label{splitcon}\end{gather}and the analogous decomposition for $M^*$ and $d-i$ (instead of $M$ and $i+1$). Then (\ref{spliteins}) holds. Conversely, if (\ref{spliteins}) holds then $F_H^{i+1}\cap F_{con}^{d-i}=0$.
\end{pro}

{\sc Proof:} In general we have $$F_{con}^{d-i}=\bi[H^d({X}_{\Gamma,alg},t_{\le i}({\mathcal {M}}\otimes_K\Omega^{\bullet}_{{X}_{\Gamma}})^{alg})\to H^d({X}_{\Gamma,alg},({\mathcal {M}}\otimes_K\Omega^{\bullet}_{{X}_{\Gamma}})^{alg})].$$Let ${\mathfrak X}_{\Gamma,alg}$ denote the ${\mathcal O}_K$-scheme (constructed in \cite{mus}) of which  ${\mathfrak X}_{\Gamma}$ is the $\pi$-adic formal completion and ${X}_{\Gamma,alg}$ the generic fibre. If $t:{\mathfrak X}_{\Gamma}\to{\mathfrak X}_{\Gamma,alg}$ and $j:{X}_{\Gamma,alg}\to{\mathfrak X}_{\Gamma,alg}$ denote the natural maps then we have a canonical transformation$${\mathbb R}j_*t_{\le i}({\mathcal {M}}^{\Gamma}\otimes_K\Omega^{\bullet}_{{X}_{\Gamma}})^{alg}\longrightarrow t_*\tau_i({\mathcal {M}}^{\Gamma}\otimes_K\Omega^{\bullet}_{{X}_{\Gamma}}).$$From (\ref{charga}) it then follows that $F_{con}^{d-i}\subset F_{\Gamma}^{d-i}$. Therefore our hypothesis implies (\ref{sumhalb}) and we may conclude as in the proof of \ref{genfro}.\hfill$\Box$\\

{\it Remarks:} (1) Observe that \ref{conj} formulates a purely algebraic approach to the splitting conjecture. In particular it invites trying to find a non-$p$-adic proof of the splitting conjecture. This remark may in particular be relevant in cases where ${X}_{\Gamma,alg}$ is the base change to $K$ of a Shimura variety defined over a global number field, see \cite{rz}. In these cases the complexes ${\mathcal {M}}^{\Gamma}\otimes_K\Omega^{\bullet}_{{X}_{\Gamma}}$ occur in the de Rham complexes of the relative de Rham cohomology (with Gauss-Manin connection) of powers of the universal abelian scheme. Using the criterion \ref{conj} one may hope to prove the splitting conjecture with global methods ! \\
(2) From the point of view of $p$-adic Hodge theory the relevance comes from the following fact: in \cite{hk} it is shown that if ${\mathcal M}$ is endowed with a structure of isoclinic $F$-isocrystal, then $H^{d}(X_{\Gamma},{\mathcal {M}}^{\Gamma}\otimes_K\Omega^{\bullet}_{{X}_{\Gamma}}    )$ receives a Frobenius structure and $F_{\Gamma}^{\bullet}$ is its corresponding (renumbered) slope filtration.

\section{The reduced Hodge spectral sequence}

\label{rehose}

For general $M$ the Hodge spectral sequence (\ref{hodss}) does not degenerate in $E_1$. For {\it rational} representations $M$ Schneider constructs a new ('reduced') Hodge spectral sequence computing $H^{*}(X_{\Gamma},{\mathcal {M}}^{\Gamma}\otimes_K\Omega^{\bullet}_{{X}_{\Gamma}}    )$ which he conjectures to degenerate in $E_1$. We discuss his conjecture in this section. 

If $\Xi_0,\ldots,\Xi_{d}$ denote the standard projective coordinate functions on $\mathbb{P}^{d}_K$, then $z_j=\Xi_j/\Xi_0$ for $j=1,\ldots,d$ are holomorphic functions on $X$. Let$$\overline{u}(z)=
\left(\begin{array}{cc}1&-z_1\quad\cdots\quad-z_d\\{0}&I_d\end{array}\right)\in{\rm SL}\sb {d+1}({\mathcal O}_X(X)).$$Let now $M$ be an irreducible $K$-rational representation of ${\rm GL}\sb {d+1}$. Suppose it has highest weight $(\lambda_0\ge\lambda_1\ge\ldots\ge\lambda_d)$. By this we mean that there exists a non zero vector $m\in M$ such that $K.m$ is stable under upper triangular matrices and generates $M$ as a $G$-representation, and such that $gm=\prod_{i=0}^da_i^{\lambda_i}m$ for all diagonal matrices $g=e_0(a_0)\cdots e_d(a_d)\in G$. Assume $\lambda_d=0$. We grade $M$ by setting$$\gr^r M=\{m\in M\quad|\quad e_0(a_0)m=a_0^{\lambda_0-r}m\,\,\mbox{for all}\,\,a_0\in K\}$$for $r\in\mathbb{Z}$, and we filter $M$ by setting$$f^rM=\bigoplus_{r'\ge r}\gr^{r'}M.$$Then $f^{\lambda_0+1}M=0$ and $f^0M=M$. We get a corresponding filtration of the constant sheaf $\underline{M}$ on $\mathfrak{X}$ and on $X$. We filter $\underline{M}\otimes_K\Omega^{j}_X$ by setting$$f^r(\underline{M}\otimes_K\Omega^{j}_X)={\mathcal O}_X.\overline{u}(z)(f^r\underline{M})\otimes_{{\mathcal O}_X}\Omega^j_X.$$We let\begin{gather}{\mathcal F}^{r,\bullet}=[f^{r}(\underline{M}\otimes_K\Omega^0_X)\longrightarrow f^{r-1}(\underline{M}\otimes_K\Omega^1_X)\longrightarrow f^{r-2}(\underline{M}\otimes_K\Omega^2_X)\longrightarrow\ldots].\label{filder}\end{gather}

It follows from \cite{schn} that this is a ${\rm SL}\sb {d+1}(K)$-stable filtration of $\underline{M}\otimes_K\Omega^{\bullet}_X$ by subcomplexes (notations and normalizations in loc. cit. are different, but equivalent). We obtain the spectral sequence\begin{gather}E_1^{r,s}=h^{r+s}({\mathcal F}^{r,\bullet}/{\mathcal F}^{r+1,\bullet})\Rightarrow h^{r+s}(\underline{M}\otimes_K\Omega^{\bullet}_X)\label{ogusfil}.\end{gather}The following is \cite{schn} Lemma 9, section 3 (observe that $X$ is a Stein space).

\begin{pro}\label{peterco} (Schneider) The terms $\underline{D}^j(M)=E_1^{\lambda_0-\lambda_j+j,\lambda_j-\lambda_0}$ for $0\le j\le d$ are the only non vanishing $E_1$-terms in (\ref{ogusfil}).\hfill$\Box$\\
\end{pro}

We define ${\rm SL}\sb {d+1}(K)$-invariant subobjects $B^j$ and $Z^j$ of $\underline{M}\otimes_K\Omega^{j}_X$ by requiring$$f^{\lambda_0-\lambda_j+1}(\underline{M}\otimes_K\Omega^{j}_X)\subset B^j\subset Z^j\subset f^{\lambda_0-\lambda_j}(\underline{M}\otimes_K\Omega^{j}_X),$$$$Z^j/f^{\lambda_0-\lambda_j+1}(\underline{M}\otimes_K\Omega^{j}_X)=\ker(\delta^j_{\lambda_j-j}),\quad\quad B^j/f^{\lambda_0-\lambda_j+1}(\underline{M}\otimes_K\Omega^{j}_X)=\bi(\delta^{j-1}_{\lambda_j-j})$$where $\delta^j_t:{\mathcal F}^{\lambda_0-t,j}/{\mathcal F}^{\lambda_0-t+1,j}\to{\mathcal F}^{\lambda_0-t,j+1}/{\mathcal F}^{\lambda_0-t+1,j+1}$ is the differential. Now \ref{peterco} implies (compare the proof of \cite{schn} Theorem 3, section 3) that\begin{gather}Z^0\longrightarrow Z^1+dB^0\longrightarrow Z^2+dB^1\longrightarrow\ldots\longrightarrow Z^d+dB^{d-1}\label{intco}\end{gather}is a subcomplex of $\underline{M}\otimes_K\Omega^{j}_X$ such that the inclusion into $\underline{M}\otimes_K\Omega^{j}_X$ is a quasiisomorphism. Moreover it implies that for any $j$ the map$$Z^j+dB^{j-1}\longrightarrow Z^j/B^j=\underline{D}^j(M)$$$$z+db\mapsto z\mod B^j$$is well defined and that if we take via these maps the quotient complex$$\underline{D}^0(M)\longrightarrow \underline{D}^1(M)\longrightarrow \underline{D}^2(M)\longrightarrow\ldots \longrightarrow \underline{D}^d(M)$$ of (\ref{intco}), then this quotient map is a quasiisomorphism, too. Hence an ${\rm SL}\sb {d+1}(K)$-equivariant isomorphism between $\underline{M}\otimes_K\Omega^{\bullet}_X$ and $\underline{D}^{\bullet}(M)$ in the derived category $D({\mathfrak{X}})$ of abelian sheaves on ${\mathfrak{X}}$. 

Let again $\Gamma<{\rm SL}\sb {d+1}(K)$ be as before. Consider the spectral sequences\begin{gather}E_1^{st}=H^{s+t}(X_{\Gamma},({\mathcal F}^{s,\bullet}/{\mathcal F}^{s+1,\bullet})^{\Gamma})\Rightarrow H^{s+t}(X_{\Gamma},{\mathcal {M}}^{\Gamma}\otimes_K\Omega^{\bullet}_{{X}_{\Gamma}} )\label{tilho}\end{gather}
\begin{gather}E_1^{st}=H^t(X_{\Gamma},\underline{D}^{s}(M)^{\Gamma})\Rightarrow H^{s+t}(X_{\Gamma},\underline{D}^{\bullet}(M)^{\Gamma})=H^{s+t}(X_{\Gamma},{\mathcal {M}}^{\Gamma}\otimes_K\Omega^{\bullet}_{{X}_{\Gamma}} ).\label{redss}\end{gather}The latter is called the 'reduced' Hodge spectral sequence computing our object of interest $H^{*}(X_{\Gamma},{\mathcal {M}}^{\Gamma}\otimes_K\Omega^{\bullet}_{{X}_{\Gamma}} )$. Let $$H^{d}(X_{\Gamma},{\mathcal {M}}^{\Gamma}\otimes_K\Omega^{\bullet}_{{X}_{\Gamma}} )={F}^0_I\supset{F}^1_{I}\supset\ldots\supset{F}^{\lambda_0+d}_{I}\supset{F}^{\lambda_0+d+1}_{I}=(0)$$be the filtration induced by $(\ref{tilho})$, let $$H^{d}(X_{\Gamma},{\mathcal {M}}^{\Gamma}\otimes_K\Omega^{\bullet}_{{X}_{\Gamma}} )=F^0_{red}\supset F^1_{red}\supset\ldots\supset F^d_{red}\supset F^{d+1}_{red}=(0)$$be the filtration induced by $(\ref{redss})$. These filtrations have the dame jumps; namely, from \ref{peterco} it follows that for all $d\ge j\ge 1$ we have\begin{gather}F_{red}^j={F}_I^{\lambda_0-\lambda_{j-1}+j}={F}_I^{\lambda_0-\lambda_{j-1}+j+1}=\ldots={F}_I^{\lambda_0-\lambda_j+j}.\label{redhod}\end{gather}

The irreducible $K$-rational ${\rm GL}\sb {d+1}$-representation$${{{}}}{M}^*=\ho_K(M,K)\otimes\det{}^{\lambda_0}$$has highest weight $({{{}}}{\lambda}_0^*\ge\ldots\ge{{{}}}{\lambda}_d^*)$ with ${{{}}}{\lambda}_d^*=0$, where ${{{}}}{\lambda}_i^*=\lambda_0-\lambda_{d-i}$ for $0\le i\le d$. A straightforward computation shows that the filtration $(f^r{{{}}}{M}^*)_r$ of ${{{}}}{M}^*$ is dual to the filtration $(f^r{M})_r$ of ${M}$, in the sense that the canonical perfect pairing$$M\times {{{}}}{M}^*\longrightarrow K$$induces perfect pairings$$\gr^{\lambda_0-j}M\times \gr^j {{{}}}{M}^*\longrightarrow K$$for any $0\le j\le \lambda_0={{{}}}{\lambda}_0^*$. These are not ${\rm SL}\sb {d+1}(K)$-equivariant objects. However, applying the ${\rm SL}\sb {d+1}(K)$-equivariance of the pairings$$\underline{M}\otimes_K\Omega^i_{X}\times{{{}}}{\underline{M}}^*\otimes_K{\Omega}^{d-i}_{X}\longrightarrow{\Omega}^d_{X},$$$$(m\otimes \eta,m^*\otimes\omega)\mapsto m^*(m)\eta\wedge\omega$$to the action of the element $\overline{u}(z)$ one deduces perfect pairings$$\frac{f^{\lambda_0-j}(\underline{M}\otimes_K\Omega^i_X)}{f^{\lambda_0-j+1}(\underline{M}\otimes_K\Omega^{i}_X)}\times \frac{f^j({{{}}}{\underline{M}}^*\otimes_K\Omega^{d-i}_X)}{f^{j+1}({{{}}}{\underline{M}}^*\otimes_K\Omega_X^{d-i})}\longrightarrow\Omega^d_X.$$Clearly they are compatible with the differential when $i$ varies, hence ${\rm SL}\sb {d+1}(K)$-equivariant perfect pairings$$\underline{D}^i(M)\times \underline{D}^{d-i}({{{}}}{M}^*)\longrightarrow \Omega^d_X.$$Passing to $\Gamma$-invariant sheaves on ${\mathfrak X}_{\Gamma}$ resp. ${X}_{\Gamma}$ we get the perfect pairing$$\underline{D}^i(M)^{\Gamma}\times \underline{D}^{d-i}({{{}}}{M}^*)^{\Gamma}\longrightarrow \Omega^d_{{X}_{\Gamma}}.$$In particular, Serre duality on the smooth projective $K$-scheme $X_{\Gamma}$ gives us perfect pairings\begin{gather}H^s(X_{\Gamma},\underline{D}^i(M)^{\Gamma})\times H^{d-s}(X_{\Gamma},\underline{D}^{d-i}({{{}}}{M}^*)^{\Gamma})\longrightarrow K.\label{serred}\end{gather} 
{\bf Conjecture:} (Schneider) For all $0\le i\le d-1$ we have \begin{gather}H^{d}(X_{\Gamma},{\mathcal {M}}^{\Gamma}\otimes_K\Omega^{\bullet}_{{X}_{\Gamma}} )=F_{red}^{i+1}\bigoplus F_{\Gamma}^{d-i}.\notag\end{gather}

\begin{satz}\label{intelde} If $\nabla(M)$ and $\nabla(M^*)$ are surjective then $F_{red}^{\bullet}=F_{H}^{\bullet}$ and\begin{gather}H^{d}(X_{\Gamma},{\mathcal {M}}^{\Gamma}\otimes_K\Omega^{\bullet}_{{X}_{\Gamma}} )=F_{red}^{i+1}\bigoplus F_{\Gamma}^{d-i}\quad\quad(0\le i\le d-1).\label{splitzwei}\end{gather}\end{satz}

{\sc Proof:} (i) We first claim that there exists a ${\rm SL}\sb {d+1}(K)$-equivariant morphism of sheaf complexes $$\nu:{\mathbb{L}}^{\bullet}(M)\longrightarrow \underline{D}^{\bullet}(M)$$ which in $D({\mathfrak{X}})$ coincides with the inclusion of sheaf complexes ${\mathbb{L}}^{\bullet}(M)\to \underline{M}\otimes_K\Omega_{{X}}^{\bullet}$, via the previous isomorphism between $\underline{M}\otimes_K\Omega^{\bullet}_X$ and $\underline{D}^{\bullet}(M)$ in $D({\mathfrak{X}})$.\\
For any $j$ denote by $d^j:\underline{M}\otimes_K\Omega_X^j\to\underline{M}\otimes_K\Omega_X^{j+1}$ the differential. By \ref{peterco} we know that $d^{i-1}$ induces a surjection\begin{gather}\underline{M}\otimes\Omega^{i-1}_X\longrightarrow\frac{\ke(d^i)}{f^{\lambda_0-\lambda_i}(\underline{M}\otimes\Omega^{i}_X)\cap\ke(d^i)}.\label{pecoco}\end{gather}Now let $\omega\in{\mathbb{L}}^{i}(M)$. Choose an element $\alpha\in \underline{M}\otimes\Omega^{i-1}_X$ which maps under (\ref{pecoco}) to the class represented by $\omega$. Then $d^{i-1}(\alpha)-\omega$ lies in $f^{\lambda_0-\lambda_i}(\underline{M}\otimes\Omega^{i}_X)\cap\ke(d^i)$ and we define $\nu(\omega)$ as its class in$$\frac{f^{\lambda_0-\lambda_i}(\underline{M}\otimes\Omega^{i}_X)\cap\ke(d^i)}{f^{\lambda_0-\lambda_i+1}(\underline{M}\otimes\Omega^{i}_X)+d^{i-1}(f^{\lambda_0-\lambda_i+1}(\underline{M}\otimes\Omega^{i-1}_X))}\subset\underline{D}^i(M).$$That $\nu$ has the stated property follows from \ref{peterco}.

(ii) Next we claim \begin{gather}H^{d}(X_{\Gamma},{\mathcal {M}}^{\Gamma}\otimes_K\Omega^{\bullet}_{{X}_{\Gamma}} )=F_{red}^{i+1}+F_{\Gamma}^{d-i}.\label{sumall}\end{gather}The map $\nu$ from (i) induces a surjective map$$H^{d}(X_{\Gamma},{\mathbb{L}}^{\bullet}(M)^{\Gamma})\longrightarrow H^{d}(X_{\Gamma},\underline{D}^{\bullet}(M)^{\Gamma})=H^{d}(X_{\Gamma},{\mathcal {M}}^{\Gamma}\otimes_K\Omega^{\bullet}_{{X}_{\Gamma}} ).$$This follows from \ref{nabkri} and the stated property of $\nu$. Let$$F_{\gamma}^{d-i}=\bi[H^{d}(X_{\Gamma},t_{\le i}{\mathbb{L}}^{\bullet}(M)^{\Gamma})\longrightarrow H^{d}(X_{\Gamma},{\mathcal {M}}^{\Gamma}\otimes_K\Omega^{\bullet}_{{X}_{\Gamma}} )],$$$$F_{\mathbb{L}}^{i+1}=\bi[H^{d}(X_{\Gamma},{\mathbb{L}}^{\bullet}(M)_{\ge i+1}^{\Gamma})\longrightarrow H^{d}(X_{\Gamma},{\mathcal {M}}^{\Gamma}\otimes_K\Omega^{\bullet}_{{X}_{\Gamma}} )].$$Then $F_{\mathbb{L}}^{i+1}\subset F_{red}^{i+1}$, again by (i), and $F_{\gamma}^{i+1}\subset F_{\Gamma}^{d-i}$, by \ref{acycanw} (since $t_{\le i}{\mathbb{L}}^{\bullet}(M)\subset \tau_i(\underline{M}\otimes_K\Omega_{{X}}^{\bullet})$). Since ${\mathbb{L}}^{\bullet}(M)=t_{\le i}{\mathbb{L}}^{\bullet}(M)\oplus{\mathbb{L}}^{\bullet}(M)_{\ge i+1}$ we get (\ref{sumall}).

(iii) (The remaining arguments are copied from the proof of \cite{iovspi} Theorem 5.4.) Let us denote by $\check{F}_{\Gamma}^{\bullet}$ and $\check{F}_{red}^{\bullet}$ the filtrations on $H^{d}(X_{\Gamma},{\mathcal{M}}^{*,{\Gamma}}\otimes_K\Omega^{\bullet}_{X_{\Gamma}})=H^{d}(X_{\Gamma},\underline{D}^{\bullet}({{{}}}{M}^*)^{\Gamma})$ analogous to the filtrations ${F}_{\Gamma}^{\bullet}$ and $F_{red}^{\bullet}$ on $H^{d}(X_{\Gamma},{\mathcal {M}}^{\Gamma}\otimes_K\Omega^{\bullet}_{{X}_{\Gamma}} )$. Here we claim\begin{align}\dim_K(H^{d}(X_{\Gamma},{\mathcal{M}}^{*,{\Gamma}}\otimes_K\Omega^{\bullet}_{X_{\Gamma}}))&=
\dim_K(H^{d}(X_{\Gamma},{\mathcal {M}}^{{\Gamma}}\otimes_K\Omega^{\bullet}_{{X}_{\Gamma}} ))\notag \\{} &=\dim_K(\check{F}_{red}^{d-i})+\dim_K(F_{red}^{i+1}).\notag\end{align}From the perfect pairings (\ref{serred}) we get perfect pairings$$H^{d}(X_{\Gamma},\underline{D}^{\bullet}(M)_{\ge i+1}^{\Gamma})\times H^{d}(X_{\Gamma},\underline{D}^{\bullet}({{{}}}{M}^*)^{\Gamma}_{\le d-i-1})\longrightarrow K$$$$H^{d}(X_{\Gamma},\underline{D}^{\bullet}(M)^{\Gamma})\times H^{d}(X_{\Gamma},\underline{D}^{\bullet}({{{}}}{M}^*)^{\Gamma})\longrightarrow K$$which commute with each other in the obvious sense. Thus $(F_{red}^{i+1})^{\bot}=\check{F}_{red}^{d-i}$ and claim (iii) follows.

(iv) The theorem is well known in case $d=1$, thus we assume $d\ge2$. From formula (\ref{pecomp}) we get$$\dim_K(F_{\Gamma}^{d-i})+\dim_K(F_{\Gamma}^{i+1})=\dim_K(H^{d}(X_{\Gamma},{\mathcal {M}}^{\Gamma}\otimes_K\Omega^{\bullet}_{{X}_{\Gamma}} )).$$This formula together with (\ref{sumall}) implies\begin{gather}\dim_K(F_{red}^{i+1})\ge\dim_K(H^{d}(X_{\Gamma},{\mathcal {M}}^{\Gamma}\otimes_K\Omega^{\bullet}_{{X}_{\Gamma}} ))-\dim_K(F_{\Gamma}^{d-i})=\dim_K(F^{i+1}_{\Gamma}).\label{ung}\end{gather}We compute\begin{align}\dim_K(F_{red}^{i+1})&=\dim_K(H^{d}(X_{\Gamma},{\mathcal{M}}^{*,{\Gamma}}\otimes_K\Omega^{\bullet}_{X_{\Gamma}}))-\dim_K(\check{F}_{red}^{d-i})\notag\\{}&\le\dim_K(\check{F}^{i+1}_{\Gamma})=\dim_K(F^{i+1}_{\Gamma}).\notag\end{align}Here the first equality follows from claim (iii), the inequality uses formula (\ref{ung}) for ${{{}}}{M}^*$ instead of $M$, and the last equality is a consequence of the formulae (\ref{pecomp}) and (\ref{dimdu}).  Altogether we see that in (\ref{ung}) we even have equality, which concludes the proof of (\ref{splitzwei}) in view of (\ref{sumall}).

(v) We have $F_{\mathbb{L}}^{i+1}\subset F_{red}^{i+1}$ and $F_{\gamma}^{d-i}\subset F_{\Gamma}^{d-i}$ (see (ii)) as well as $F_{\mathbb{L}}^{i+1}\subset F_{H}^{i+1}$. On the other hand $F_{\gamma}^{d-i}+F_{\mathbb{L}}^{i+1}=H^{d}(X_{\Gamma},{\mathcal {M}}^{\Gamma}\otimes_K\Omega^{\bullet}_{{X}_{\Gamma}})$ by the surjectivity of $\nabla(M)$. Since we have $F_{red}^{i+1}\cap F_{\Gamma}^{d-i}=0=F_{H}^{i+1}\cap F_{\Gamma}^{d-i}$ we find $F_{\gamma}^{d-i}=F_{\Gamma}^{d-i}$ and $F_{red}^{d-i}=F_{\mathbb{L}}^{i+1}=F_{H}^{d-i}$.\hfill$\Box$\\

Denote by ${\mathbb L}_D^{\bullet}(M)$ the subsheaf complex of $\underline{D}^{\bullet}(M)$ on ${\mathfrak X}\otimes k$ defined by $${\mathbb L}_D^{i}(M)(U)=\ke[\underline{D}^{i}(M)(]U[)\longrightarrow\underline{D}^{i+1}(M)(]U[)]$$for open $U\subset{\mathfrak X}\otimes k$. The inclusion ${\mathbb L}_D^{\bullet}(M)^{\Gamma}\to\underline{D}^{\bullet}(M)^{\Gamma}$ induces a map$$\theta(M):H^d(\mathfrak{X}_{\Gamma},{\mathbb L}_D^{\bullet}(M)^{\Gamma})\longrightarrow H^d(\mathfrak{X}_{\Gamma},\underline{D}^{\bullet}(M)^{\Gamma}).$$

\begin{satz}\label{redcri} (a) If $\theta(M)$ and $\theta(M^*)$ are surjective then we have the decomposition (\ref{splitzwei}).\\(b) The following two statements (i) and (ii) are equivalent:\\(i) For any $i,j$ the following map is bijective:$$H^j(\mathfrak{X}_{\Gamma},{\mathbb L}_D^{i}(M)^{\Gamma})\longrightarrow H^j(\mathfrak{X}_{\Gamma},\underline{D}^{i}(M)^{\Gamma})$$(ii) We have (\ref{splitzwei}), and the reduced Hodge spectral sequence $(\ref{redss})$ degenerates in $E_1$.  
\end{satz}

{\sc Proof:} (a) Proposition \ref{acycanw} also holds if $\tau_i(\underline{M}\otimes\Omega_X^{\bullet})$ is replaced by$$\tau_i\underline{D}^{\bullet}(M)=[\underline{D}^{0}(M)\longrightarrow\ldots\longrightarrow \underline{D}^{i-1}(M)\longrightarrow{\mathbb L}_D^{i}(M)\longrightarrow\ldots].$$Indeed, in view of the quasiisomorphism of {\it sheaf} complexes $\underline{D}^{\bullet}(M)\cong\underline{M}\otimes\Omega_X^{\bullet}$ this version is in fact reduced to \ref{acycanw}. As in \ref{genfro} we therefore obtain\begin{gather}F_{\Gamma}^{d-i}=\bi[H^d({\mathfrak{X}}_{\Gamma},\tau_i\underline{D}^{\bullet}(M)^{\Gamma})\to H^{d}(\mathfrak{X}_{\Gamma},\underline{D}^{\bullet}(M)^{\Gamma})]\label{gamred}\end{gather}and we get claim (a) just as in \ref{nabkri} and/or \ref{intelde}. Claim (b) is then also clear, again using (\ref{gamred}).\hfill$\Box$\\

{\it Remark:} In \cite{schn} it is conjectured that $(\ref{redss})$ always degenerates in $E_1$. Thus \ref{redcri} (b)(i) should be a sufficient {\it and necessary} condition to prove the decomposition (\ref{splitzwei}) !  

\begin{kor}\label{gammafi} Suppose that $M$ is a $K$-rational $G$-representation with small weights.\\(a) The splitting in Corollary \ref{hosssp} is given by the filtration $F_{\Gamma}^{\bullet}$:\begin{gather}H^{d}(X_{\Gamma},{\mathcal {M}}^{\Gamma}\otimes_K\Omega^{\bullet}_{{X}_{\Gamma}}    )=F_H^{i+1}\bigoplus F_{\Gamma}^{d-i}\quad\quad(0\le i\le d-1).\label{smaspli}\end{gather}(b) We have $F_H^{\bullet}=F_{red}^{\bullet}$ in $H^{d}(X_{\Gamma},{\mathcal {M}}^{\Gamma}\otimes_K\Omega^{\bullet}_{{X}_{\Gamma}})$.
\end{kor}

{\sc Proof:} This is the combination of Corollary \ref{nabkri} and Theorem \ref{intelde} with Corollary \ref{hosssp}. We may pass to the base field extension $K\to\dot{K}$. Then we have inclusions of sheaf complexes$${\mathcal Log}^{\bullet}({\mathcal M}^0_{{\mathcal O}_{\dot{K}}})^{\Gamma}\otimes_{{\mathcal O}_{\dot{K}}} \dot{K}\longrightarrow{\mathbb{L}}^{\bullet}(M)^{\Gamma}\otimes_K\dot{K}\longrightarrow{\mathcal {M}}\otimes_K\Omega^{\bullet}_{{X}_{\Gamma}}\otimes_K\dot{K}$$and similarly for $M^*$ instead of $M$. Thus Corollary \ref{hosssp} implies that $\nabla(M)$ and $\nabla(M^*)$ are surjective and (a) follows from \ref{nabkri} and (b) follows from \ref{intelde}.\hfill$\Box$\\

{\it Remarks:} (1) The decomposition (\ref{smaspli}) was proven for the trivial representation $M=K$ for the first time by Iovita and Spiess \cite{iovspi}. Our present proof appears to provide a geometric underpinning of the one given in \cite{iovspi}.\\(2) The degeneration of the Hodge spectral sequence (\ref{hodss}) is of course well known for $M=K$ and $\kara(K)=0$. On the other hand, for more general $K$-rational representations $M$ than those considered in \ref{hosssp} it can not be expected to degenerate (see \cite{schn}).\\(3) Let $I$ be a $K[\Gamma]$-module (with $\dim_KI<\infty$) which contains a $\Gamma$-stable free ${\mathcal O}_K$-lattice $I^0$. Let $\sigma\in{\rm Gal}(K/{\mathbb Q}_p)$, let $M$ be a $K$-rational $G$-representation and let $M_{\sigma}$ denote $M$ but with the $K$-vector space structure twisted by $\sigma$ --- then $G$ acts on $M_{\sigma}$ again by $K$-linear automorphisms. Everything we did in this paper with the local system defined by $M$ carries over to the local system defined by $I\otimes_K M_{\sigma}$: simply replace every occurence of ${\mathcal M}^0_{{\mathcal O}_{\dot{K}}}$ by ${I}^0\otimes_{{\mathcal O}_K}{\mathcal M}^0_{{\mathcal O}_{\dot{K}},\sigma}$.\\(4) Let $\breve{K}$ denote the (completed) maximal unramified extension of $K$. The formal scheme $\mathfrak{X}\times\spf({\mathcal O}_{\breve{K}})$ carries a certain universal $G$-equivariant formal group ${\mathcal G}$, see \cite{rz}. If $K={\mathbb Q}_p$ the de Rham complex ${\mathfrak E}\otimes\Omega_{X\dot{\otimes}\breve{K}}^{\bullet}$ of its relative Dieudonne module ${\mathfrak E}$ (as a filtered convergent $F$-isocrystal on $\mathfrak{X}\times\spf({\mathcal O}_{\breve{K}})$) can be identified with a sum of $d+1$ copies of $K^{d+1}\otimes_K\Omega_{X\dot{\otimes}_K\breve{K}}^{\bullet}$, filtered as in (\ref{filder}) and with isotypical Frobenius action of slope $d/(d+1)$. From our results in \cite{hk} it follows that the filtration $F_{\Gamma}^{\bullet}$ on $H^d(X_{\Gamma}\otimes \breve{K},{\mathfrak E}^{\Gamma}\otimes\Omega_{X_{\Gamma}\otimes \breve{K}}^{\bullet})$ is the (renumbered) slope filtration. Hence \ref{gammafi} states that the slope filtration on $H^d(X_{\Gamma}\otimes \breve{K},{\mathfrak E}^{\Gamma}\otimes\Omega_{X_{\Gamma}\otimes \breve{K}}^{\bullet})$ is opposite to the Hodge filtration. By the comparison isomorphisms of $p$-adic Hodge theory this is a statement on the cohomology of the relative Tate module of the $\Gamma$-quotient of ${\mathcal G}$.


\end{document}